\documentclass[11pt]{article}
\usepackage{amsfonts}
\usepackage{amssymb}
\usepackage{amsthm}
\usepackage{amsmath}
\usepackage{latexsym}
\usepackage[T1]{fontenc}
\usepackage[latin1]{inputenc}
\usepackage{graphicx}
\usepackage{color}

\def \bop {\noindent\textbf{Proof. }}
\def \eop {\hbox{}\nobreak\hfill \vrule width 2mm height 2mm depth 0mm}

\newtheorem{definition}{Definition}[section]
\newtheorem{theorem}{Theorem}[section]
\newtheorem{proposition}{Proposition}[section]
\newtheorem{lemma}{Lemma}[section]
\newtheorem{remark}{Remark}[section]
\newtheorem{corollary}{Corollary}[section]

\def \eop {\hbox{}\nobreak\hfill \vrule width 2.0mm height 1.8mm depth 0mm
\par \goodbreak \smallskip}
\input{tcilatex}
\begin{document}
\author{\textsf{K. Bahlali}$^{\mathrm{a}}$\textsf{, M. Eddahbi}$^{\mathrm{b}%
} $ and \textsf{Y. Ouknine}$^{\mathrm{c}}$ \\
$^{\mathrm{a}}${\small \textit{Université de Toulon, IMATH, EA 2134, 83957
La Garde, France}} \\
{\small \ e--mail: bahlali@univ-tln.fr} \\
$^{\mathrm{b}}${\small \textit{Université Cadi Ayyad, Faculté des Sciences
et Techniques, }}\\
{\small \textit{Département de Mathématiques, B.P. 549, Marrakech, Maroc.}}\\
{\small \ e--mail: m.eddahbi@uca.ma} \\
$^{\mathrm{c}}${\small \textit{Université Cadi Ayyad, Faculté des Sciences
Semlalia,}}\\
{\small \ \textit{Département de Mathématiques, B.P. 2390\ Marrakech, Maroc.}
}\\
{\small \ e--mail: ouknine@uca.ma}}
\title{Quadratic BSDEs with $\mathbb{L}^2$--terminal data Existence results, Krylov's estimate and It\^o--Krylov's formula \thanks{%
Partially supported by Marie Curie Initial Training Network (ITN)
\textquotedblleft Deterministic and Stochastic Controlled Systems
and Applications\textquotedblright\ (PITN-GA-2008-213841-2), PHC
Volubilis MA/10/224 and PHC Tassili 13MDU887}}

\maketitle

\begin{abstract}
In a first step, we establish the existence (and sometimes the uniqueness)
of solutions for a large class of quadratic backward stochastic differential
equations (QBSDEs) with continuous generator and a merely square integrable
terminal condition. Our approach is different from those existing in the
literature. Although we are focused on QBSDEs, our existence result also
covers the BSDEs with linear growth, keeping $\xi$ square integrable in both
cases. As byproduct, the existence of viscosity solutions is established for
a class of quadratic partial differential equations (QPDEs) with a square
integrable terminal datum. In a second step, we consider QBSDEs with
measurable generator for which we establish a Krylov's type a priori
estimate for the solutions. We then deduce an It\^o--Krylov's change of
variable formula. This allows us to establish various existence and
uniqueness results for classes of QBSDEs with square integrable terminal
condition and sometimes a merely measurable generator. Our results show, in
particular, that neither the existence of exponential moments of the
terminal datum nor the continuity of the generator are necessary to the
existence and/or uniqueness of solutions for quadratic BSDEs.
Some comparison theorems are also established for solutions of a class of
QBSDEs.
\end{abstract}



\section{Introduction}

Let $(W_{t})_{0\leq t\leq T}$ be a $d$--dimensional Brownian motion defined
on a complete probability space $(\Omega ,\mathcal{F},\mathbb{P})$. We
denote by $(\mathcal{F}_{t})_{0\leq t\leq T})$ the natural filtration of $W$
augmented with $\mathbb{P}$--negligible sets. Let $H(t,\omega ,y,z)$ be a
real valued $\mathcal{F}_{t}$--progressively measurable process defined on $%
[0,\ T]\times \Omega \times \mathbb{R}\times \mathbb{R}^{d}$. Let $\xi $ be
an $\mathcal{F}_{T}$--measurable $\mathbb{R}$--valued random variable.
In this paper, we consider a one dimensional BSDE of the form,
\begin{equation}
Y_{t}=\xi +\int_{t}^{T}H(s,Y_{s},Z_{s})ds-\int_{t}^{T}Z_{s}dW_{s},\ \ 0\leq
t\leq T  \tag{$eq(\xi ,H)$}
\end{equation}%
The data $\xi $ and $H$ are respectively called the terminal condition and
the coefficient or the generator of the BSDE $eq(\xi ,H)$.

A BSDE is called quadratic if its generator has at most a quadratic growth
in the $z$ variable. \newline
For given real numbers $a$ and $b$, we set \ $a\wedge b:=\min (a,b)$, \ $%
a\vee b:=\max (a,b)$, \ $a^{-}:=\max (0,-a)$ \ and \ $a^{+}:=\max (0,a)$. We
also define,\\
$\mathcal{W} _{1,\,loc}^{2}$ := the Sobolev space of (classes) of
functions $u$ defined on $\mathbb{R}$ such that both $u$ and its generalized
derivatives $u^{\prime }$ and $u^{\prime \prime }$ belong to $\mathbb{L}%
_{loc}^{1}(\mathbb{R})$.\\
$\mathcal{S}^{2}$ := the set of continuous, $\mathcal{F}_{t}$%
--adapted processes $\varphi$ such that $$\mathbb{E}\sup\limits_{0\leq
t\leq T}|\varphi_{t}|^{2} < \infty. $$
$\mathcal{M}^{2}$ := the space of $\mathcal{F}_{t}$--adapted processes $
\varphi$ satisfying $\mathbb{E} \int_{0}^{T}|\varphi_{s}|^{2}ds<+\infty. $\\
$\mathcal{L}^{2}$ := the space of $\mathcal{F}_{t}$ --adapted processes $%
\varphi$ satisfying $\int_{0}^{T}|\varphi_{s}|^{2}ds<+\infty \ \mathbb{P}%
\text{--a.s.} \label{inq2} $

\begin{definition}
\label{SolutionL2} A solution to BSDE $eq(\xi ,H)$ is an $\mathcal{F}_{t} $%
--adapted processes $(Y,Z)$ which satisfy the BSDE $eq(\xi ,H)$ for each $%
t\in [ 0,T]$ and such that $Y$ is continuous \ and \ $\int_0^T\vert
Z_{s}\vert^2ds < \infty$ \ $\mathbb{P}$--a.s., that is $(Y,\, Z) \in
\mathcal{C}\times\mathcal{L}^{2}$, where $\mathcal{C}$ is the space of
continuous processes.
\end{definition}

The first results on the existence of solutions to QBSDEs were obtained
independently in \cite{Kob} and in \cite{DHO} by two different methods. The
approach developed in \cite{Kob} is based on the monotone stability of
QBSDEs and consists to find bounded solutions. Later, many authors have
extended the result of \cite{Kob} in many directions, see $e.\,g.$ \cite{BE,
BH1, EH1, LSM2, Tevz}. For instance, in \cite{BH1}, the existence of
solutions was proved for QBSDEs in the case where the exponential moments of
the terminal datum are finite. In \cite{Tevz}, a fixed point method is used
to directly show the existence and uniqueness of a bounded solution for
QBSDEs with a bounded terminal datum and a (so--called) Lipschitz--quadratic
generator. More recently, a monotone stability result for quadratic
semimartingales was established in \cite{BE} then applied to derive the
existence of solutions to QBSDEs in the framework of exponential
integrability of the terminal data. The generalized stochastic QBSDEs were
studied in \cite{EH1} under more or less similar assumptions on the terminal
datum. Applications of QBSDEs in financial mathematics are also given in
\cite{BE} with a large bibliography in this subject.

It should be noted that all the previous papers in QBSDEs were developed in
the framework of continuous generators and bounded terminal data or at least
having finite exponential moments. It is natural to ask the following questions :

1) Are there quadratic BSDEs that have solutions without assuming the
existence of exponential moments of the terminal datum ? If yes, in what space these solutions lie ?

2) Are there quadratic BSDEs with measurable generator that have solutions
without assuming the existence of exponential moments of the terminal datum ?
If yes, in what space these solutions lie ?

\vskip0.1cm The present paper gives positive answers to these questions. It
is a development and a continuation of our announced results \cite{BEO1}. We
do not aim to generalize the previous papers on QBSDEs, but our goal is to give
another point of view (on solving QBSDEs) which allows us to establish the
existence of solutions, in the space $\mathcal{S}^{2}\times \mathcal{M}^{2}$,
for a large class of QBSDEs with a square integrable terminal datum. Next, in
order to deal with QBSDEs with measurable generator, we had to establish a
Krylov's type a priori estimate and an It\^o--Krylov's formula for the
solutions of general QBSDEs.

\vskip 0.15cm To begin, let us give a simple example which is covered by the present paper but, to the best of our knowledge, is not covered by the previous results.
This example shows that the existence of exponential moments of the terminal
datum is not necessary to the unique solvability of BSDEs in $\mathcal{S}%
^{2}\times \mathcal{M}^{2}$. Assume that,

\vskip 0.2cm\noindent\textbf{(H1)} \ \ \ \ \ \textit{$\xi$ is square
integrable.}

\vskip 0.2cm\noindent Let $f: \mathbb{R}\longmapsto \mathbb{R}$ be a given
continuous function with compact support, and set $M:= \sup_{y\in \mathbb{R}%
}|f(y)|$. The BSDE $eq(\xi ,f(y)|z|^2)$ is then of quadratic growth since $%
|f(y)|z|^2|\leq M|z|^{2}$.
Let $u(x):=\int_{0}^{x}\exp \left( 2\int_{0}^{y}f(t)dt\right) dy $. If $(Y,Z)
$ is a solution to the BSDE $eq(\xi ,f(y)|z|^2)$, then It\^o's formula applied
to $u(Y_{t})$ shows that,
\begin{equation*}
u(Y_{t})=u(\xi )-\int_{t}^{T}u^{\prime }(Y_{s})Z_{s}dW_{s}
\end{equation*}%
If we set $\bar{Y_{t}}:=u(Y_{t})$ and $\bar{Z_{t}}:=u^{\prime }(Y_{t})Z_{t}$%
, then $(\bar{Y}, \bar{Z})$ solves the BSDE
\begin{equation*}
\bar{Y_{t}}=u(\xi )-\int_{t}^{T}\bar{Z_{s}}dW_{s}
\end{equation*}%
Since both $u$ and its inverse are $\mathcal{C}^{2}$ smooth functions which
are globally Lipschitz and one to one from $\mathbb{R}$ onto $\mathbb{R}$,
we then deduce that the BSDE $eq(\xi ,f(y)|z|^2)$ admits a solution (resp. a
unique solution) if and only if the BSDE $eq(u(\xi ),0)$ admits a solution
(resp. a unique solution). The BSDE $eq(u(\xi ),0)$ has a unique solution in
$\mathcal{S}^{2}\times \mathcal{M}^{2}$ whenever $u(\xi) $ is merely square
integrable. According to the properties of $u$ and its inverse, $%
u(\xi)$ is square integrable if and only if $\xi$ square integrable.
Therefore, even when all the exponential moments are infinite the QBSDE $%
eq(\xi ,f(y)|z|^2)$ has a unique solution which lies in $\mathcal{S}%
^{2}\times \mathcal{M}^{2}$. Note that, since the $sign$ of $f$ is not
constant, our example also shows that the convexity of the generator is not
necessary to the uniqueness. Assume now that $\xi$ is merely $\mathcal{F}_T$
measurable, but not necessarily integrable. According to Dudley's
representation theorem \cite{Dud}, one can show as previously (by using the
above transformation $u$) that when $f$ is continuous and with compact
support, the BSDE $eq(\xi ,f(y)|z|^2)$ has at least one solution $(Y,Z)$
which belongs to $\mathcal{C}\times\mathcal{L}^{2}$.

\vskip 0.2cm In the first part of this paper, we establish the existence of
solutions for a large class of QBSDEs having a continuous generator and a
merely square integrable terminal datum. The generator $H$ will satisfy

\hskip 3cm $|H(s,y,z)| \leq (a+b|y|+c|z|+f(|y|)|z|^{2})$

\noindent where $f$ is some continuous and globally integrable function on $\mathbb{R}$
(hence can not be a constant) and $a, b, c$ are some positive constants.

\noindent Our approach consists to deduce the solvability of a BSDE (without
barriers) from that of a suitable QBSDE with two Reflecting barriers whose
solvability is ensured by \cite{EH1}. This allows us to control the
integrability we impose to the terminal datum. In other words, this idea
can be summarized as follows:  When $|H(s,y,z)| \leq
(a+b|y|+c|z|+f(|y|)|z|^{2})$, the existence of solutions for the QBSDE $eq(\xi, H)$ can be
deduced from the existence of solutions to the
 QBSDE driven by the dominating generator $a+b|y|+c|z|+f(|y|)|z|^{2}$. Using the
transformation $u$ (defined in the above first example), we show that the
solvability of the QBSDE $eq(\xi, \ a+b|y|+c|z|+f(|y|)|z|^{2})$ is
equivalent to the solvability of a BSDE without quadratic term which is more
easily solvable. We also prove that the uniqueness of solutions holds for
the class of QBSDEs $eq(\xi ,f(y)|z|^2)$ under the ${\mathbb{L}}^2$--integrability condition
on the terminal data. It is worth to notice that the
existence results of \cite{BE, BH1, Kob, LSM1, LSM2} can be obtained by our
method.  We mention that, in contrast to the most previous papers on QBSDEs,
our result also cover the BSDEs with linear growth (by putting $f=0$). It
therefore provides a unified treatment for quadratic BSDEs and those of
linear growth, keeping $\xi$ square integrable in both cases.


\vskip 0.15cm

In the second part of this paper, we begin by proving the Krylov inequality
for the solutions of general QBSDEs from which we deduce the It\^o--Krylov
formula, $i.e.$ we show that the It\^o change of variable formula holds for
$u(Y_t)$ whenever $Y$ is a solution of a QBSDE, $u$ is of class $\mathcal{C}%
^1$ and the second generalized derivative of $u$ merely belongs to $\mathbb{L%
}_{loc}^{1}( \mathbb{R})$.

We then use this change of variable formula to establish the existence (and
sometimes the uniqueness) of solutions. To explain more precisely how we get
our second aim, let us consider the following assumption,

\vskip0.2cm\noindent \textbf{(H2)} \ \textit{There exist a positive
stochastic process $\eta \in \mathbb{L}^{1}([0,\ T]\times \Omega )$ and a
locally integrable}

\hskip 0.5cm \textit{function $f$ such that for every $(t,\omega ,y,z)$,}
\begin{equation*}
|H(t,y,z)|\leq \eta _{t}+|f(y)||z|^{2} \ \ \ \ \ \mathbb{P}\otimes dt \ a.e.
\end{equation*}

We first use the occupation time formula to show that if assumption \textbf{%
(H2)} holds, then for any solutions $(Y,Z)$ of the BSDE $eq(\xi,H)$, the
time spend by $Y$ in a Lebesgue negligible set is negligible with respect to
the measure $|Z_{t}|^{2}dt$. That is, the following Krylov's type estimate
holds for any positive measurable function $\psi $,
\begin{equation}
\mathbb{E}\int_{0}^{T\wedge \tau _{R}}\psi (Y_{s})|Z_{s}|^{2}ds\leq
C\left\Vert \psi \right\Vert _{\mathbb{L}^{1}([-R,R])},  \label{kryintro}
\end{equation}
where $\tau _{R}$ is the first exit time of $Y$ from the interval $[-R,R]$
and $C$ is a constant depending on $T$, $||\xi ||_{\mathbb{L}^{1}(\Omega )}$
and $||f||_{\mathbb{L}^{1}\left( [-R,R]\right) }$.

\noindent We then deduce (by assuming \textbf{(H1)}--\textbf{(H2)}) that :
if $(Y,Z)$ is a solution to the BSDE $eq(\xi,H)$ which belongs to $\mathcal{S}
^{2}\times \mathcal{L}^{2}$, then for any function $\varphi$ $\in$ $\mathcal{%
C}^{1}(\mathbb{R})\cap \mathcal{W} _{1,\,loc}^{2}(\mathbb{R})$ the following
change of variable formula holds true,
\begin{equation}  \label{ikintro}
\varphi(Y_{t})=\varphi(Y_{0})+\int_{0}^{t}\varphi^{\prime} (Y_{s})dY_{s}+%
\frac{1}{2} \int_{0}^{t}\varphi^{\prime\prime} (Y_{s})|Z_{s}|^{2}ds
\end{equation}

Inequality (\ref{kryintro}) as well as formula (\ref{ikintro}) are
interesting in their own and can have potential applications in BSDEs. They
are established here with minimal conditions on the data $\xi$ and $H$.
Indeed, it will be shown that these formulas hold for QBSDEs with a merely
measurable generator. For formula (\ref{ikintro}) we require that the
terminal datum is square integrable, while for inequality (\ref{kryintro})
we do not need any integrability condition on the terminal datum. Notice that,
although the inequality (\ref{kryintro}) can be established by adapting the
method developed by Krylov, which is based on partial differential equations
\cite{Krlivre80} (see also \cite{Bflows, BM, Mel, Kr87}), the proof we give
here is purely probabilistic and more simple.

\vskip 0.15cm As application, we establish the existence of solutions in $%
\mathcal{S} ^{2}\times \mathcal{M}^{2}$ for the classes of QBSDEs  $eq(\xi
,\, f(y)|z|^{2})$ and $eq(\xi ,\, a+by+cz+f(y)|z|^{2})$ assuming merely that
$f$ is globally integrable and $\xi$ is square integrable. Remark that, when $f
$ is not continuous, the function $u(x) :=\int_{0}^{x}\exp \left(
2\int_{0}^{y}f(t)dt\right) dy$ is not of $\mathcal{\ C}^2$--class and the
classical It\^o's formula can not be applied. Nevertheless, when $f$ belongs
to $\mathbb{L}^{1}(\mathbb{R})$, the function $u$ belongs to the space $%
\mathcal{C }^{1}(\mathbb{R})\cap \mathcal{W}_{1,\,loc}^{2}(\mathbb{R})$ and
hence formula (\ref{ikintro}) can be applied to $u$. Our strategy consists
then to use the idea we developed in the first part to show the existence of
a minimal and a maximal solution for BSDE $eq(\xi ,\,
a+b|y|+c|z|+f(|y|)|z|^{2})$.

A comparison theorem is also proved for two BSDEs of type $eq(\xi
,f(y)|z|^{2})$ whenever we can compare their terminal data and $a.e.$ their
generators. We then deduce the uniqueness of solutions for the BSDEs $eq(\xi
, f(|y|)|z|^{2})$ when $\xi$ is square integrable and $f$ belongs to $%
\mathbb{\mathbb{L}}^1(\mathbb{R})$. That is, even when $f$ is defined merely
$a.e.$, the uniqueness holds. This gives a positive answer to question 2. In
particular, the QBSDE $eq(\xi , H)$ has a unique solution $(Y,Z)$ which
belongs to $\mathcal{S}^{2}\times \mathcal{M}^{2}$ when $\xi $ is merely
square integrable and $H$ is one of the following generators:

\vskip
0.15cm \ $H_1(y,z) := \sin(y)|z|^2$ if $y\in[-\pi,\,\frac{\pi}{2}]$ \ and $\
H_1(y,z) := 0 $ \ otherwise,

\vskip 0.15cm \ $H_2(y,z) : =( \mathbf{1}_{[a,b]}(y)- \mathbf{1}%
_{[c,d]}(y))|z|^2$ for a given $a<b$ and $c<d$,

\vskip 0.15cm \ $H_3(y,z) := \frac{1}{(1+y^{2})\sqrt{|y|}}|z|^2$ if $y \neq 0
$ \ and $\ H_2(y,z) := 1 $ \ otherwise.

\vskip 0.15cm \noindent It should be noted that the generator $H_3(y,z)$ is
neither continuous nor locally bounded and the QBSDE $eq(\xi , H_3)$ has a
unique solution in $\mathcal{S}^{2}\times \mathcal{M}^{2}$ when $\xi $ is
merely square integrable.

\vskip 0.15cm We finally consider the BSDE $eq(\xi ,H)$. We assume that $\xi$
is square integrable and $H$ is continuous in $(y,z)$ and $|H(s, y, z)| \leq a+b|y|+c|z|+f(|y|)|z|^{2})$, with $f$ merely globally integrable and
locally bounded but not necessarily continuous. Although one can argue as in
the first part to obtain the existence of solutions from the solutions of a
suitable Reflected BSDE, we give a different proof which is based on a
classical comparison theorem and an appropriate localization by a suitable
dominating process which is derived from the extremal solutions of the QBSDE
$eq(\xi ,(a+b|y|+c|z|+f(|y|)|z|^{2}))$. This allows us to construct a
suitable sequence of BSDEs $eq(\xi_n ,H_n)$ whose localized (i.e. stopped)
solutions converge to a solution of the BSDE $eq(\xi ,H)$.

\vskip 0.15cm In the third part, we establish the existence of viscosity
solutions for a class of non--divergence form semilinear PDEs with quadratic
nonlinearity in the gradient variable. This is done with a continuous generator and an unbounded terminal
datum. It surprisingly turns out that there is
a gap between the BSDEs and the classical formulation of their associated
semilinear PDEs (see Remark \ref{gap}, section 5). Observe that the class of quadratic PDEs we study in
this paper can be used as a simplified model in some incomplete financial
markets, see e.g. \cite{DJ}.

\vskip 0.15cm The paper is organized as follows. In section 2, we study the
QBSDEs with a continuous generator and a square integrable terminal datum.
Krylov's estimate and It\^o --Krylov's formula for QBSDEs are established in
section 3. The solvability of a class of QBSDEs with measurable generator is
studied in section 4. An application to the existence of viscosity solutions
for Quadratic PDEs associated to the Markovian QBSDE $eq(\xi ,f(y)|z|^{2})$
is given in section 5.

%




\section{QBSDEs with $\mathbb{L}^2$ terminal data and continuous generators}

\vskip 0.15cm We will establish the solvability in $\mathcal{S}^{2}\times
\mathcal{M}^{2}$ for some BSDEs with a square integrable terminal data and a
continuous generator. Our method consists to construct a solution of a BSDE
without barriers from a solution of a suitable BSDE with two Reflecting
barriers. More precisely : Assuming that $\xi$ is square integrable and $f$
is continuous and globally integrable on $\mathbb{R}$, we first establish
the existence of a minimal and a maximal solution in $\mathcal{S}^{2}\times
\mathcal{M}^{2}$ for the BSDE $eq(\xi ,\, a+b|y|+c|z|+f(|y|)|z|^{2})$. An
next, we consider the BSDE $eq(\xi ,H)$ with \ $|H(s,y,z)| \leq
(a+b|y|+c|z|+f(|y|)|z|^{2})$.  We then use the minimal solution of BSDE $%
eq(-\xi^- ,-(a+b|y|+c|z|+f(|y|)|z|^{2}))$ and the maximal solution of BSDE $%
eq(\xi^+ ,a+b|y|+c|z|+f(|y|)|z|^{2})$ as barriers, and apply the result of
\cite{EH1} to get the existence of a solution which stays between these two
barriers. We finally deduce the solvability of $eq(\xi ,H)$ by proving that
the increasing stochastic processes, which force the solutions to stay
between the barriers, are equal to zero.


\vskip 0.15cm The following lemma is needed for the sequel of the paper. It
allows us to eliminate the additive quadratic term.

\begin{lemma}
\label{zvonkin0C} Let $f$ be continuous and belongs to $\mathbb{L}^{1}(%
\mathbb{R})$. The function
\begin{equation}
u(x):=\int_{0}^{x}\exp \left( 2\int_{0}^{y}f(t)dt\right) dy
\end{equation}
has the following properties,

(i) \ $u\in \mathcal{C}^{2}(\mathbb{R})$ \ and satisfies the equation \ $%
\frac{1}{2}u^{\prime \prime }(x)-f(x)u^{\prime }(x)=0$, in $\mathbb{R}$.

(ii) \ $u$ is a one to one function from $\mathbb{R}$ onto $\mathbb{R}$.

(iii) The inverse function $u^{-1}$ belongs to $\mathcal{C}^{2}(\mathbb{R})$.

(iv) \ $u$ is a quasi--isometry, that is there exist two positive constants $%
m $ and $M$ such that,

\hskip 1.5cm for any \ $x,y\in \mathbb{R}$, \ $m\left\vert x-y\right\vert
\leq \left\vert u(x)-u(y)\right\vert \leq M\left\vert x-y\right\vert $
\end{lemma}


\subsection{The equation $eq(\protect\xi, f(y)|z|^2)$}

\begin{remark}
Let $\xi $ be an $\mathcal{F}_{T}$--measurable random variable. According to
Dudley \cite{Dud}, there exists a (non necessary unique) $\mathcal{F}_{t}$
--adapted process $(Z_{t})_{0\leq t\leq T}$ such that $%
\int_{0}^{T}|Z_{s}|^{2}ds<\infty $ $\mathbb{P}$--$a.s$ and $\xi
=\int_{0}^{T}Z_{s}dW_{s}$. The process $(Y_{t})_{0\leq t\leq T}$ defined by $%
Y_{t}=\xi -\int_{0}^{t}Z_{s}dW_{s}$ is $\mathcal{F}_{t}$--adapted and
satisfies the equation $eq(\xi ,0)$. This solution $(Y_{t},Z_{t})_{0\leq
t\leq T}$ is not unique. However, if we assume $\xi \in L^{2}(\Omega )$ then
the solution $(Y,Z)$ is unique and $Y_{t}=\mathbb{E}\left[ \xi /\mathcal{F}%
_{t}\right] $.
\end{remark}

\vskip 0.1cm The following proposition shows that the exponential moment of $%
\xi$ is not needed to obtain the existence and uniqueness of the solution to
quadratic BSDEs.

\begin{proposition}
\label{f(y)|z|^2C} (i) Assume \textbf{\ (H1)} be satisfied. Let $f$ be a
continuous and integrable function. Then the BSDE $eq(\xi,f(y)|z|^{2})$ has
a unique solution in $\mathcal{S}^{2}\times \mathcal{L}^{2}$ (resp. in $%
\mathcal{S}^{2}\times \mathcal{M}^{2}$) if and only if the BSDE $%
eq(u(\xi),0) $ has a unique solution in $\mathcal{S}^{2}\times \mathcal{L}%
^{2}$ (resp. in $\mathcal{S}^{2}\times \mathcal{M}^{2}$).

(ii) In particular, the BSDE $eq(\xi,f(y)|z|^{2})$ has a unique solution $%
(Y, \, Z)$ which belongs to $\mathcal{S}^2\times\mathcal{M}^2$.
\end{proposition}

\bop Let $u$ be the function defined in Lemma \ref{zvonkin0}. Theorem \ref%
{Ito-Krylov} and Lemma \ref{zvonkin0} allow us to show that, $(Y_{t},Z_{t})$
is the unique solution of BSDE $eq(\xi ,\,f(y)|z|^2)$ if and only if $(\bar{Y%
},\,\bar{Z}):=(u(Y_{t}),\,u^{\prime }(Y_{t})Z_{t})$ is the unique solution
to BSDE $eq(u(\xi ),0)$. We shall prove assertion $(ii)$. Since $\xi $ is
square integrable then $u(\xi )$ is square integrable too. Therefore $%
eq(u(\xi ),\,0)$ has a unique solution in $\mathcal{S}^2\times\mathcal{M}^2$%
. Assertion $(ii)$ follows now from assertion $(i)$. \eop


\subsection{The equation $eq(\protect\xi, a+b|y|+c|z|+f(y)|z|^{2})$}

The BSDE under consideration in this subsection is,
\begin{equation}
Y_{t}=\xi
+\int_{t}^{T}(a+b|Y_{s}|+c|Z_{s}|+f(Y_{s})|Z_{s}|^{2})ds-%
\int_{t}^{T}Z_{s}dW_{s}  \label{edsrab1}
\end{equation}%
where $a, \, \, b,\, \, c \in \mathbb{R}$ and $f:\mathbb{R}\longmapsto
\mathbb{R}$.

\vskip0.15cm We refer to BSDE (\ref{edsrab1}) as equation $eq(\xi
,a+b|y|+c|z|+f(y)|z|^{2})$.

\begin{proposition}
\label{equivabfC} Assume that \textbf{\ (H1)} holds. Assume also
that $f$ is continuous and globally integrable on $\mathbb{R}$. Let
$u$ be the function defined in Lemma \ref{zvonkin0}. Then the BSDE
$eq(\xi ,\,a+b|y|+ c| z|+f(y)|z|^{2})$ has at least one
solution. Moreover all
solutions of $eq(\xi ,\,a+b|y|+ c| z|+f(y)|z|^{2})$ are in $\mathcal{S}%
^{2}\times\mathcal{M}^{2}$.
\end{proposition}


\bop It\^o's formula applied to the function $u$ (which is defined
in Lemma \ref{zvonkin0C}) shows that $(Y_t,Z_t)$ is 
solution to the BSDE $eq(\xi ,a+b|y|+ c| z|+f(y)|z|^2)$ if and only
if $(\bar Y_t, \bar Z_t):=(u(Y_t),u^{\prime }(Y_t)Z_t)$ is a
 solution to the BSDE
$eq(u(\xi ),(a+b|u^{-1}(\bar y)|) u^{\prime }[u^{-1}(\bar y)]+c|\bar
z|)$. We shall prove the existence of solutions to BSDE $eq(u(\xi
),(a+b|u^{-1}(\bar y)|) u^{\prime }[u^{-1}(\bar y)]+c|\bar z|)$. By Theorem \ref{Ito-Krylov} we have
\begin{equation}  \label{edsrLSM}
\bar Y_{t}=\bar \xi +\int_{t}^{T}G(\bar Y_{s}, \bar Z_{s})ds-
\int_{t}^{T}\bar Z_{s}dW_{s}
\end{equation}
where $G(\bar y,\bar z):=(a+b|u^{-1}(\bar y)|)u^{\prime }[u^{-1}(\bar y)] +
c|\bar z|$.

\noindent From Lemma \ref{zvonkin0C}, we deduce that the generator $G$ is
continuous and with linear growth, and the terminal condition $\bar\xi
:=u(\xi )$ is square integrable (since Assumption \textbf{(H1)}). Hence,
according to Lepeltier \& San-Martin \cite{LSM1}, the BSDE (\ref{edsrLSM})
has at least one solution in $\mathcal{S}^{2}\times\mathcal{M}^{2}$. To
complete the proof, it is enough to observe that the function $u$ defined in
Lemma \ref{zvonkin0C} is strictly increasing. \eop


\vskip 0.3cm\noindent \textbf{Alternative proof to Proposition \ref%
{equivabfC}.} \\ In the previous proof of Proposition
\ref{equivabfC}, we had to use the Lepeltier \& San-Martin result
\cite{LSM1} in order to quickly deduce the existence of
solutions which belong to $\mathcal{S}^{2}\times\mathcal{M}^{2}$
whenever the terminal datum $\xi$ is square integrable. This fact
will be proved below by using an alternative proof which is in
adequacy with the spirit of the present paper. To this end, we use a
result on two barriers Reflected QBSDEs obtained by Essaky $\&$
Hassani in \cite{EH2} which establishes the existence of solutions
for reflected QBSDEs without assuming any integrability condition
on the terminal datum. For the self--contained, we state the result of \cite%
{EH2} in the following theorem.


\begin{theorem}
\label{EH} (\cite{EH2}, Theorem 3.2). Let $L$ and $U$ be continuous
processes and $\xi$ be a $\mathcal{F}_T$ measurable random variable. Assume
that

1) \hskip 0.2cm for every $t\in[0, \ T]$, \ $L_t \leq U_{t}$

2) \hskip 0.2cm $L_T\leq\xi\leq U_T$.

3) \hskip 0.2cm there exists a continuous semimartingale which passes
between the barriers $L$ and $U$.

4) \hskip 0.2cm $H$ is continuous in $(y,z)$ and satisfies for every $%
(s,\omega)$, every $y\in [L_s(\omega), U_s(\omega)]$ and every $z\in \mathbb{%
R}^d$.

\hskip 2cm $|f(s, \omega, y, z )| \leq \eta_s(\omega)+ C_s(\omega)|z|^2 $

\vskip 0.2cm\noindent where $\eta \in \mathbb{L}^{1}([0, T]\times\Omega)$
and $C$ is a continuous process.

\vskip 0.2cm\noindent Then, the following RBSDE has a minimal and a maximal
solution.
\begin{equation}  \label{eq000}
\left\{
\begin{array}{ll}
(i) & Y_{t}=\xi + \displaystyle\int_{t}^{T}H(s,Y_{s},Z_{s})ds  - \displaystyle%
\int_{t}^{T}Z_{s}dB_{s} \\
& + \displaystyle%
\int_{t}^{T}dK_{s}^+ - \displaystyle\int_{t}^{T}dK_{s}^- \; \mbox{for all }\; t\leq T \\
(ii) & \forall t\leq T,\,\, L_t \leq Y_{t}\leq U_{t},\quad \\
(iii) & \displaystyle\int_{t}^{T}( Y_{t}-L_{t}) dK_{t}^+ = \displaystyle%
\int_{t}^{T}( U_{t}-Y_{t}) dK_{t}^-=0,\,\, \mbox{a.s.}, \\
(iv) & K_0^+ =K_0^- =0, \,\,\,\, K^+, K^- \,\,%
\mbox{are continuous
nondecreasing}. \\
(v) & dK^+ \bot dK^-%
\end{array}
\right.
\end{equation}
\end{theorem}
\vskip 0.2cm We are now in the position to give our alternative
proof to Proposition \ref{equivabfC}. \newline Note that, since $u$
is strictly increasing, we then only need to prove the existence of
a minimal and a maximal solutions for the BSDE (\ref{edsrLSM}).
Since $\xi$ is square integrable, then according to Lemma
\ref{zvonkin0C} the terminal condition $\bar\xi :=u(\xi )$ is also
square integrable. Once again, by using Lemma \ref{zvonkin0C}, one
can show that the generator $G$ of the BSDE (\ref{edsrLSM}) is
continuous and with linear growth. Indeed, since $u(0) = 0$ and
$u^{\prime }$ is bounded by $M$ (Lemma \ref{zvonkin0C} $ (iv)$), we
have
\begin{align}  \label{borneGbar}
G(y,z)& =(a+b|u^{-1}(y)|)u^{\prime }[u^{-1}( y)] + c| z|  \notag \\
& \leq M a + m M b|y| + c| z| :=g(y,z)
\end{align}
where $m$ and $M$ are the constants which appear in assertion $(iv)$
of Lemma \ref{zvonkin0C}. \newline Since the function $g(y, z) := M
a + m M b|y| + c| z|$ is uniformly Lipschitz and with linear growth
in $(y,z)$, then according to Pardoux $\&$ Peng result \cite{PP1},
the BSDEs $eq(-\bar\xi^-, -g)$ and $eq(\bar\xi^+, g)$ have unique
solutions in $\mathcal{S}^{2}\times\mathcal{M}^{2}$, which we
respectively denote by $(Y^{-g},Z^{-g})$) and $(Y^{g},Z^{g})$. Note
that $ Y^{-g} $ is negative and $Y^{g}$ is positive. Using Theorem
\ref{EH} (with $L =Y^{-g}$, $U=Y^{g}$, $\eta_t = M a + m M
b(|Y_t^{-g}| + |Y_t^{g}|) + c^2$, and $C_t = 1 $), we deduce
that the Reflected BSDE 
\begin{equation}  \label{RBSDEHdominé}
\left\{
\begin{array}{ll}
(i) & Y_{t}=\xi +\displaystyle\int_{t}^{T}H(s,Y_{s},Z_{s})ds
-\displaystyle
\int_{t}^{T}Z_{s}dB_{s}\,, t\leq T, \\
& +\displaystyle
\int_{t}^{T}dK_{s}^+ -\displaystyle\int_{t}^{T}dK_{s}^- \; \mbox{for all }\; t\leq T\\
(ii) & \forall \ t\leq T,\,\, Y_t^{-g} \leq Y_{t}\leq Y_t^{g}, \\
(iii) & \displaystyle\int_{0}^{T}( Y_{t}-Y_t^{-g}) dK_{t}^+ =
\displaystyle
\int_{0}^{T}( Y_t^{g}-Y_{t}) dK_{t}^-=0,\,\, \mbox{a.s.}, \\
(iv) & K_0^+ =K_0^- =0, \,\,\,\, K^+, K^- \,\,
\mbox{are continuous
nondecreasing}. \\
(v) & dK^+ \bot dK^-
\end{array}
\right.
\end{equation}
has at least one solution $ (Y, Z, K^{ +}, K^{ -})$ and $(Y,Z)$ belongs to $\mathcal{C}\times\mathcal{L}^2$.

We shall show that $dK^+ = dK^- = 0$.  Since $Y_t^{g}$ is a
solution to the BSDE $ eq(\bar\xi^+,g)$, then Tanaka's formula
applied to $(Y_t^{g}-Y_t)^+$ shows that
\begin{align*}
(Y_{t}^g-Y_{t})^{+} & =(Y_{0}^g-Y_{0})^{+}+\int_{0}^{t} \mathbf{1}
_{\{Y_{s}^g > Y_s\}} [H(s,Y_s,Z_s)-g(s,Y_s^g,Z_s^g)]ds \\
& \ + \int_{0}^{t} \mathbf{1}_{\{Y_{s}^g > Y_s\}}(dK_{s}^+-dK_{s}^-) +
\int_{0}^{t} \mathbf{1}_{\{Y_{s}^g > Y_s\}}(Z_{s}^g-Z_s)dW_s \\
& \ + L_t^0(Y^g-Y)
\end{align*}
where $L_t^0(Y^g-Y)$ is the local time at time $t$ and level $0$ of the
semimartingale $(Y^g-Y)$.

\vskip 0.15cm\noindent Since $Y^{g}\geq Y$, then $(Y_{t}^{g}-Y_{t})^{+} =
(Y_{t}^{g}-Y_{t})$. Therefore, identifying the terms of $%
(Y_{t}^{g}-Y_{t})^{+}$ with those of $(Y_{t}^{g}-Y_{t})$ and using the fact
that:
\begin{equation*}
\mathbf{1}-\mathbf{1}_{\{Y_{s}^{g}>Y_{s}\}}=\mathbf{1} _{\{Y_{s}^{g}\leq
Y_{s}\}}=\mathbf{1}_{\{Y_{s}^{g}=Y_{s}\}},
\end{equation*}
we obtain,
\begin{equation*}
(Z_s - Z_s^g)\mathbf{1}_{\{Y_{s}^{g}=Y_{s}\}} = 0 \ \text{for} \ a.e. \ (s, \omega)
\end{equation*}
Using the previous equalities, one can show that
\begin{eqnarray*}
&& \int_{0}^{t}\mathbf{1}_{\{Y_{s}^{g}=Y_{s}
\}}(dK_{s}^{+}-dK_{s}^{-})=L_t^0(Y^g-Y) \\ &&+ \int_{0}^{t}\mathbf{1}
_{\{Y_{s}^{g}=Y_{s} \}}[g(s,Y_{s}^{g},Z_{s}^{g})-H(s,Y_{s},Z_{s})]ds
\end{eqnarray*}
Since $\int_{0}^{t}\mathbf{1}_{\{Y_{s}^{g}=Y_{s}\}}dK_{s}^{+}=0$, it holds
that
\begin{eqnarray*}
0&\leq & L_t^0(Y^g-Y)+\int_{0}^{t}\mathbf{1}_{\{Y_{s}^{g}=Y_{s}
\}}[g(s,Y_{s}^{g},Z_{s}^{g})-H(s,Y_{s},Z_{s})]ds\\ &&=-\int_{0}^{t}\mathbf{1}
_{\{Y_{s}^{g}=Y_{s}\}}dK_{s}^{-}\leq 0
\end{eqnarray*}
Hence,
$\int_{0}^{t}\mathbf{1}_{\{Y_{s}^{g}=Y_{s}\}}dK_{s}^{-}=0$, which
implies that \ $dK^{-}=0$. Arguing symmetrically, one can show that
$ dK^{+}=0 $. Therefore $(Y,Z)$ is a solution to the (non reflected)
BSDE $ eq(\xi ,H)$. Moreover $Y$ belongs to $\mathcal{S} ^{2}$ since both $Y^{g}$ and $Y^{-g}$ belong to
$\mathcal{S} ^{2}$. Remember that $G$ is of linear
growth and $Y$ belongs to $\mathcal{S}^{2}$, then using standard
arguments of BSDEs, we deduce that $Z$ belongs to
$\mathcal{M}^2$. This completes the ``alternative proof of
Proposition \ref{equivabfC}". \eop

\begin{remark}
The previous proof also constitute an alternative proof to the result of
Lepeltier $\&$ San-Martin for the existence of a minimal and a maximal
solution to BSDEs with continuous and at most of linear growth generator. It
is worth noting that the idea consists to construct a solution of a BSDE
with linear growth from a solution of a Reflected Quadratic BSDE.
\end{remark}

\subsection{The BSDE $eq (\protect\xi ,H)$, with $\vert H(s,y,z)\vert \leq a
+ b\vert y\vert + c\vert z\vert + f(\vert y\vert)|z|^2$}

\noindent Consider the assumptions,

\vskip 0.2cm\noindent \textbf{(H4)} For $a.e.\ (s,\omega )$, $H$ is
continuous in $(y,z)$

\vskip 0.2cm\noindent \textbf{(H5)} There exist positive real numbers $a,b,
c $ such that for every $s,y,z$ \vskip 0.2cm \hskip 2cm $|H(s,y,z)|\leq
a+b|y|+ c\vert z\vert + f(|y|)|z|^{2} \ := \ g(y,z)$, \vskip 0.2cm\noindent
where $f$ is some positive continuous and integrable function.

\begin{theorem}
\label{BSDEHdominéC} Assume that (\textbf{H1}), (\textbf{H4}) and (\textbf{H5%
}) are fulfilled. Then, the BSDE $eq(\xi ,H)$ has at least one solution \ $%
(Y,Z)$ in $\mathcal{S}^{2}\times\mathcal{M}^{2}$.
\end{theorem}


\vskip 0.15cm\noindent \textbf{Proof of Theorem \ref{BSDEHdominéC}}. The
idea is close to the above ``Alternative proof of proposition \ref{equivabfC}%
" and consists to derive the existence of solution for the BSDE without
reflection from solutions of a suitable $2$--barriers Reflected BSDE. Put \ $%
g(y,z) := a+b|y|+ c|z|+ f(|y|)|z|^{2}$. According to Proposition \ref%
{equivabfC}, let $(Y^{g},Z^{g})$ be a solution of BSDE
$eq(\xi^{+} ,g)$ and $(Y^{-g},Z^{-g})$) be a solution
of BSDE $eq(-\xi^{-} ,-g) $. We know by Proposition \ref{equivabfC}
that $(Y^{g}$ and $Y^{-g}$)
belong to $\mathcal{S}^{2}$. Using Theorem \ref{EH} (with $L =Y^{-g}$, $%
U=Y^{g}$, $\eta_t = a + b(|Y_t^{-g}| + |Y_t^{g}|) + c^2$, and $C_t = 1 + %
\sup_{s\leq t}\sup_{\alpha\in [0,1]}|f(\alpha
Y_s^{-g} +(1-\alpha) Y_s^{g})|$), we deduce the existence of solution $%
(Y,Z,K^+,K^-)$ to the following Reflected BSDE, such that $(Y,Z)$ belongs to
$\mathcal{C}\times\mathcal{L} ^{2}$.

\begin{equation}  \label{RBSDEHdominé2}
\left\{
\begin{array}{ll}
(i) & \!\!Y_{t}=\xi\! +\!\displaystyle\int_{t}^{T}H(s,Y_{s},Z_{s})ds  \!-\!\displaystyle %
\int_{t}^{T}Z_{s}dB_{s}, \\ & \quad \!+\!\displaystyle %
\int_{t}^{T}dK_{s}^+\! -\!\displaystyle\int_{t}^{T}dK_{s}^- \;\; \mbox{for all }\; t \leq T\\
(ii) & \forall \ t\leq T,\, Y_t^{-g} \leq Y_{t}\leq Y_t^{g}, \\
(iii) & \displaystyle\int_{0}^{T}( Y_{t}-Y_t^{-g}) dK_{t}^+ \!=\!  \displaystyle %
\int_{0}^{T}( Y_t^{g}-Y_{t}) dK_{t}^-=0,\,\, \mbox{a.s.}, \\
(iv) & K_0^+ =K_0^- =0, \,\, K^+, K^- \,\,
\mbox{are continuous
nondecreasing}. \\
(v) & dK^+ \bot dK^-%
\end{array}
\right.
\end{equation}
Arguing as in the proof of Proposition \ref{equivabfC} we end-up with $dK^+ = dK^- = 0$.

Therefore $(Y,Z)$ satisfies the (non reflected) BSDE $eq(\xi ,H)$%
. Note that since both $Y^{g}$ and $Y^{-g}$ belong to $\mathcal{S}^{2}$,
then $Y\in \mathcal{S}^{2}$ belongs to $\mathcal{S}^{2}$ too.

In order to complete the proof of Theorem \ref{BSDEHdominéC}, it remains to
show that $Z$ belongs to $\mathcal{M}^2$. To this end, we need the following
lemma. 

\begin{lemma}
\label{zvonkin1C} Let $f$ be continuous and integrable function on $\mathbb{R}$. Set
$$K(y):=\int _{0}^{y} \exp\left(-2\int_{0}^{x}f(r)dr\right)dx.$$ The function $$
u(x):=\int_{0}^{x}K(y) \exp\left(2\int_{0}^{y}f(t)dt\right)dy$$ satisfies following
properties:

\vskip 0.15cm \noindent $(i)$ $u$ belongs to $\mathcal{C}^{2}(\mathbb{R})$, and, $%
u(x)\geq0$ and $u^{\prime }(x)\geq0$ for $x\geq0$.

\vskip 0.15cm \noindent Moreover, $u$ satisfies, for a.e. $x$, \ $\frac{1}{2}%
u^{\prime\prime}(x)- f(x)u^{\prime}(x)= \frac{1}{2}$.

\vskip 0.15cm \noindent $(ii)$ The map $x\longmapsto v(x) := u(\vert x\vert)$ belongs
to $\mathcal{C} ^{2}(\mathbb{R})$, and $v^{\prime}(0) = 0$.

\vskip 0.15cm \noindent $(iii)$ There exist $c_1>0,\ c_2>0$ such that for every $x\in%
\mathbb{R}$, $u(\vert x\vert) \leq c_1 \vert x\vert^2$ and $u^{\prime}(\vert
x\vert) \leq c_2\vert x\vert $.
\end{lemma}

\noindent We now prove that $Z$ belongs to $\mathcal{M}^2$.

For $N>0$, let $\tau _{N}:=\inf \{t>0:|Y_{t}|+\int_{0}^{t}|v^{\prime
}(Y_{s})|^{2}|Z_{s}|^{2}ds\geq N\}\wedge T$. Set $\mbox{sgn}(x)=1$ if $x\geq
0 $ and $\mbox{sgn}(x)=0$ if $x<0$. Let $u$ be the function defined in Lemma %
\ref{zvonkin1C} and $v(y):=u(|y|)$. Since $v$ belongs to $\mathcal{C}^{2}(%
\mathbb{R})$, then using It\^o's formula it holds that for every $t\in [0,T]$,
\begin{align*}
u(|Y_{t\wedge \tau _{N}}|)& \!=\!u(|Y_{0}|)\!+\!\int_{0}^{t\wedge \tau _{N}}\!\!\left[
\frac{1}{2}u^{\prime \prime }(|Y_{s}|)|Z_{s}|^{2}-\mbox{sgn}(Y_{s})u^{\prime
}(|Y_{s}|)H(s,Y_{s},Z_{s})\right]\!\!ds \\
& \ \ \ \ \ +\int_{0}^{t\wedge \tau _{N}}\mbox{sgn}(Y_{s})u^{\prime
}(|Y_{s}|)Z_{s}dW_{s}
\end{align*}%
Passing to expectation and using successively assumption \textbf{(H5)} and
Lemma \ref{zvonkin1C}, we get for any $N>0$
\begin{align*}
&u(|Y_{0}|) =u(|Y_{t\wedge \tau _{N}}|)\\ &\!+\! \int_{0}^{t\wedge \tau _{N}}\!\!\left[%
\mbox{sgn}(Y_{s})u^{\prime }(|Y_{s}|)H(s,Y_{s},Z_{s})-\frac{1}{2}u^{\prime
\prime }(|Y_{s}|)|Z_{s}|^{2}\right]\!\!ds \\
&\!\! \leq \!u(|Y_{t\wedge \tau _{N}}|)\\& +\int_{0}^{t\wedge \tau _{N}}\!\!\left[%
u^{\prime }(|Y_{s}|)(a+b|Y_{s}|+c|Z_{s}|+f(|Y_{s}|)|Z_{s}|^{2})-\frac{1}{2}%
u^{\prime \prime }(|Y_{s}|)|Z_{s}|^{2}\right]\!\!ds \\
& \leq u(|Y_{t\wedge \tau _{N}}|)+\int_{0}^{t\wedge \tau _{N}}\!\!\left[%
u^{\prime }(|Y_{s}|)(a+b|Y_{s}|+c|Z_{s}|)-\frac{1}{2}|Z_{s}|^{2}\right]\!\!ds \\
&\!\! \leq \!u(|Y_{t\wedge \tau _{N}}|)+\int_{0}^{t\wedge \tau _{N}}\!\!\left[%
u^{\prime }(|Y_{s}|)(a+b|Y_{s}|)+c\,u^{\prime }(|Y_{s}|)|Z_{s}|-\frac{1}{2}%
|Z_{s}|^{2}\right]\!\!ds \\
&\!\! \leq \! u(|Y_{t\wedge \tau _{N}}|)\\&+\int_{0}^{t\wedge \tau _{N}}\!\!\left[%
u^{\prime }(|Y_{s}|)(a+b|Y_{s}|)+4[c\,u^{\prime }(|Y_{s}|)]^{2}+\frac{1}{4}%
|Z_{s}|^{2}-\frac{1}{2}|Z_{s}|^{2}\right]\!\!ds
\end{align*}%
Hence,
\begin{equation*}
\frac{1}{4}\mathbb{E}\!\int_{0}^{t\wedge \tau _{N}}\!\!|Z_{s}|^{2}ds\leq
u(|Y_{0}|)+\mathbb{E}\!\int_{0}^{T}\!\!\left[(a+b|Y_{s}|)u^{\prime
}(|Y_{s}|)+4c^{2}(u^{\prime }(|Y_{s}|))^{2}\right]\!ds
\end{equation*}%
We successively use Lemma \ref{zvonkin1C} -$(iii)$, the fact that the
process $Y$ belongs to $\mathcal{S}^{2}$ and Fatou's lemma, to show that $%
\mathbb{E}\int_{0}^{T}|Z_{s}|^{2}ds<\infty .$ Theorem \ref{BSDEHdominéC} is
proved. \eop


\section{Krylov's estimates and It\^o--Krylov's formula in QBSDEs}

\begin{remark}
$(i)$ The Krylov estimate for QBSDEs is obtained with minimal conditions.
Indeed, the generator $H$ will be assumed merely measurable and the terminal
condition $\xi$ merely integrable.

\vskip 0.1cm (ii) It is worth noting that the change of variable formula we
will establish here for the solutions of QBSDEs is valid although the
martingale part of $Y$ can degenerate. Actually, the martingale part of $Y$
can degenerate with respect to the Lebesgue measure but remains
nondegenerate with respect to the measure $|Z_{t}|^{2}dt$.

\vskip 0.2cm (iii) The Krylov estimate for QBSDE we state in the next
proposition can be established by using Krylov's method \cite{Krlivre80}
(see also \cite{Bflows, BM, Kr87, Mel}), which is based on partial
differential equations. The proof we give here is probabilistic and very
simple. It is based on the time occupation formula.
\end{remark}

\subsection{Krylov's estimates in QBSDEs.}


\begin{proposition}
\label{Krylovloc} (Local estimate) Assume \textbf{(H2)} holds. Let $(Y,Z)$
be a solution of the BSDE $eq(\xi,H)$ and assume that $\int%
\nolimits_{0}^{T}|H(s,Y_{s},Z_{s})|ds<\infty \ \ \mathbb{P}$--a.s.  Then,
there exists a positive constant $C$ depending on $T$, $R$ and $\left\Vert
f\right\Vert _{\mathbb{L}^{1}([-R,R])}$ such that for any nonnegative
measurable function $\psi$,
\begin{align*}
\mathbb{E}\int_{0}^{T\wedge\tau_{R}}\psi(Y_{s})|Z_{s}|^{2}ds\leq C\left\Vert
\psi\right\Vert _{\mathbb{L}^{1}([-R,R])},
\end{align*}
where $\tau_{R}:=\inf\{t>0:|Y_{t}|\geq R\}$.
\end{proposition}

\bop  Without loss of generality, we can and assume that $\eta =0$ in
assumption \textbf{(H2)}. Set $\tau _{N}^{\prime }:=\inf \{t>0,\ \int_{0}^{t}|Z_{s}|^{2}ds\geq N\}$, $%
\tau _{M}^{\prime \prime }:=\inf \{t>0,\
\int_{0}^{t}|H(s,Y_{s},Z_{s})|ds\geq M\}$, and put $\tau :=\tau _{R}\wedge
\tau _{N}^{\prime }\wedge \tau _{M}^{\prime \prime }$. Let $a$ be a real
number such that $a\leq R$. By Tanaka's formula, we have
\begin{eqnarray*}
(Y_{t\wedge \tau }-a)^{-}& =&(Y_{0}-a)^{-}-\int_{0}^{t\wedge \tau }\mathbf{1}%
_{\{Y_{s}<a\}}dY_{s}+\frac{1}{2}L_{t\wedge \tau }^{a}(Y) \\
& =&(Y_{0}-a)^{-}-\int_{0}^{t\wedge \tau }\mathbf{1}_{\{Y_{s}<a%
\}}H(s,Y_{s},Z_{s})ds\\&&+\int_{0}^{t\wedge \tau }\mathbf{1}_{\{Y_{s}<a%
\}}Z_{s}dW_{s}+\frac{1}{2}L_{t\wedge \tau }^{a}(Y)
\end{eqnarray*}%
Since the map $y\mapsto (y-a)^{-}$ is Lipschitz, we obtain
\begin{eqnarray}  \label{LT}
\frac{1}{2}L_{t\wedge \tau }^{a}(Y)&\leq &|Y_{t\wedge \tau
}-Y_{0}|+\int_{0}^{t\wedge \tau }\mathbf{1}_{\{Y_{s}<a\}}H(s,Y_{s},Z_{s})ds\nonumber \\ && -%
\int_{0}^{t\wedge \tau }\mathbf{1}_{\{Y_{s}<a\}}Z_{s}dW_{s}
\end{eqnarray}
Passing to expectation, we obtain
\begin{equation}
\sup_{a}\mathbb{E}\left[ L_{t\wedge \tau }^{a}(Y)\right] \leq 4R+2M
\label{estimationtempslocal}
\end{equation}
Since $-R\leq Y_{t\wedge \tau }\leq R$ for each $t$, then $%
Support(L^a_{\cdot}(Y_{\cdot \wedge \tau }))\subset [ -R,R]$. Therefore,
using inequality (\ref{LT}), assumption \textbf{(H2)} and the time
occupation formula, we get
\begin{align*}
\frac{1}{2}L_{t\wedge \tau }^{a}(Y)& \leq |Y_{t\wedge \tau
}-Y_{0}|\!+\!\!\int_{0}^{t\wedge \tau }\!\!\mathbf{1}_{\{Y_{s}<a%
\}}|f(Y_{s})||Z_{s}|^{2}ds\!-\!\!\int_{0}^{t\wedge \tau }\!\!\mathbf{1}%
_{\{Y_{s}<a\}}Z_{s}dW_{s} \\
& \leq |Y_{t\wedge \tau }-Y_{0}|\!+\!\!\int_{0}^{t\wedge \tau }\!\!\mathbf{1}%
_{\{Y_{s}<a\}}|f(Y_{s})|d\langle Y\rangle _{s}\!-\!\!\int_{0}^{t\wedge \tau
_{R}^{N}}\!\!\mathbf{1}_{\{Y_{s}<a\}}Z_{s}dW_{s} \\
& \leq |Y_{t\wedge \tau }-Y_{0}|\!+\!\!\int_{-R}^{a}|f(x)|L_{t\wedge \tau
_{R}^{N}}^{x}(Y)dx\!-\!\!\int_{0}^{t\wedge \tau }\!\!\mathbf{1}_{\{Y_{s}<a%
\}}Z_{s}dW_{s}
\end{align*}%
Passing to expectation, we obtain
\begin{equation*}
\frac{1}{2}\mathbb{E}\left[ L_{t\wedge \tau }^{a}(Y)\right] \leq \mathbb{E}%
|Y_{t\wedge \tau }-Y_{0}|+\int_{-R}^{a}|f(x)|\mathbb{E}\left[ L_{t\wedge
\tau }^{x}(Y)\right] dx<\infty
\end{equation*}%
Hence, by inequality (\ref{estimationtempslocal}) and Gronwall lemma we get
\begin{align*}
\mathbb{E}\left[ L_{t\wedge \tau }^{a}(Y)\right] & \leq 2\mathbb{E}%
(|Y_{t\wedge \tau }-Y_{0}|)\exp \left( 2\int_{-R}^{a}|f(x)|dx\right) \\
& \leq 2\mathbb{E}(|Y_{t\wedge \tau }-Y_{0}|)\exp \left( 2||f||_{\mathbb{L}%
^{1}([-R,R])}\right) \\
& \leq 4R\exp (2||f||_{\mathbb{L}^{1}([-R,R])})
\end{align*}%
Passing to the limit on $N$ and $M$ (having in mind that $\tau :=\tau
_{R}\wedge \tau _{N}^{\prime }\wedge \tau _{M}^{\prime \prime }$) and using
Beppo--Levi theorem we get
\begin{equation}
\mathbb{E}\left[ L_{t\wedge \tau _{R}}^{a}(Y)\right] \leq 4R\exp (2||f||_{%
\mathbb{L}^{1}([-R,R])})  \notag  \label{Gronwall4}
\end{equation}%
Let $\psi $ be an arbitrary positive function. We use the previous
inequality to show that
\begin{align*}
\mathbb{E}\int_{0}^{T\wedge \tau _{R}}\psi (Y_{s})|Z_{s}|^{2}ds& =\mathbb{E}%
\int_{0}^{T\wedge \tau _{R}}\psi (Y_{s})d\langle Y\rangle _{s} \\
& \leq \mathbb{E}\int_{-R}^{R}\psi (a)L_{T\wedge \tau _{R}}^{a}(Y)da \\
& \leq \int_{-R}^{R}\psi (a)\mathbb{E}L_{T\wedge \tau _{R}}^{a}(Y)da \\
& \leq 4R\exp (2||f||_{\mathbb{L}^{1}([-R,R])})\left\Vert \psi \right\Vert _{%
\mathbb{L}^{1}([-R,R])}
\end{align*}%
Proposition \ref{Krylovloc} is proved. \eop



\vskip 0.15cm We now consider the following assumption.

\vskip 0.2cm\noindent \textbf{(H3)} \textit{The function $f$, defined in
assumption \textbf{(H2)}, is globally integrable on $\mathbb{R}$.}

\vskip 0.2cm\noindent Arguing as previously, one can prove the following
\textit{global estimate}.

\begin{corollary}
\label{Krylov} Assume that \textbf{(H1)}, \textbf{(H2)} and \textbf{(H3)} are satisfied. \\ Let $(Y,Z)\in \mathcal{S}^{2}\times \mathcal{L}^{2}$ be
a solution of BSDE $eq(\xi ,H)$. Assume moreover that $\int%
\nolimits_{0}^{T}|H(s,Y_{s},Z_{s})|ds<\infty \ \ \mathbb{P}$--a.s.  Then,
there exists a positive constant $C$ depending on $T$, $\left\Vert
\xi\right\Vert _{\mathbb{L}^{1}(\Omega)}$, $\left\Vert f\right\Vert _{%
\mathbb{L}^{1}(\mathbb{R} )}$ and $\mathbb{E}(\sup_{t\leq T}\vert Y_t\vert)$
such that, for any nonnegative measurable function $\psi$,
\begin{equation}  \label{Kryoglobale}
\mathbb{E}\int_{0}^{T}\psi(Y_{s})|Z_{s}|^{2}ds\leq C\left\Vert
\psi\right\Vert _{\mathbb{L}^{1}(\mathbb{R})}
\end{equation}
In particular,
\begin{equation*}
\mathbb{E}\int_{0}^{T\wedge \tau _{R}}\psi (Y_{s})|Z_{s}|^{2}ds\leq
C\left\Vert \psi \right\Vert _{\mathbb{L}^{1}([-R,R])}.
\end{equation*}
where $\tau _{R}:=\inf \{t>0:|Y_{t}|\geq R\}$.
\end{corollary}


\subsection{An It\^o--Krylov's change of variable formula in BSDEs}

In this subsection we shall establish an It\^o--Krylov's change of variable
formula for the solutions of one dimensional BSDEs. This will allows us to
treat some QBSDEs with measurable generator. Let's give a summarized
explanation on It\^o--Krylov's formula. The It\^o change of variable formula
expresses that the image of a semimartingale, by a $\mathcal{C}^{2}$--class
function, is a semimartingale. When
\begin{equation*}
X_{t}:=X_{0}+\int_{0}^{t}\sigma (s,\omega )dW_{s}+\int_{0}^{t}b(s,\omega
)ds\
\end{equation*}%
is an It\^o's semimartingale, the so--called It\^o--Krylov's formula
(established by N.V. Krylov) expresses that if $\sigma \sigma ^{\ast }$ is
uniformly elliptic, then It\^o's formula also remains valid when $u$ belongs
to $\mathcal{W}_{p,\,loc}^{2}$ with $p$ strictly more large than the
dimension of the process $X$. Here $\mathcal{W} _{p,\,loc}^{2}$ denotes the
Sobolev space of (classes) of functions $u$ defined on $\mathbb{R}$ such
that both $u$ and its generalized derivatives $u^{\prime }$, $u^{\prime
\prime }$ belong to $L_{loc}^{p}(\mathbb{R})$. The It\^o--Krylov formula was
extended in \cite{ BM} to continuous semimartingales $%
X_{t}:=X_{0}+M_{t}+V_{t}$ with a non degenerate martingale part and some
additional conditions. The non degeneracy means that the matrix of the
increasing processes $\left<M^{i},M^{j}\right>$ is uniformly elliptic.

\begin{theorem}
\label{Ito-Krylov} Assume that (\textbf{H1}) and (\textbf{H2}) are
satisfied. Let $(Y,Z)$ be a solution of BSDE $eq(\xi,H)$ in $\mathcal{S}%
^{2}\times \mathcal{L}^{2}$. Then, for any function $u$ belonging to the
space $\mathcal{C}^{1}(\mathbb{R})\cap\mathcal{W}_{1,\, loc}^{2}(\mathbb{R }%
) $, we have
\begin{equation}  \label{ikformula}
u(Y_{t})=u(Y_{0})+\int_{0}^{t}u^{\prime} (Y_{s})dY_{s}+\frac{1}{2}
\int_{0}^{t}u^{\prime\prime} (Y_{s})|Z_{s}|^{2}ds
\end{equation}
\end{theorem}

\vskip 0.15cm Using Sobolev's embedding theorem and Lemma \ref{Ito-Krylov},
we get,

\begin{corollary}
\label{Ito-Krylov2} Assume (\textbf{H1}) and (\textbf{H2}) be satisfied. Let
$(Y,Z)$ be a solution of BSDE $eq(\xi ,H)$ in \ $\mathcal{S}^{2}\times
\mathcal{L}^{2}$. Then, for any function $u\in \mathcal{W}_{p,\, loc}^{2}(%
\mathbb{R})$ with $p>1$, the formula (\ref{ikformula}) remains valid.
\end{corollary}

\vskip 0.2cm\noindent\textbf{Proof of Theorem \ref{Ito-Krylov}.} For $R>0$,
let $\tau_{R}:=\inf\{t>0: |Y_{t}|\geq R\}$. Since $\tau_{R}$ tends to
infinity as $R$ tends to infinity, it then suffices to establish the formula
for $u(Y_{t\wedge\tau_{R}})$. Using Proposition \ref{Krylovloc} , the term $%
\int_{0} ^{t\wedge\tau_{R}}u^{\prime\prime}(Y_{s})Z_{s}^{2}ds$ is well
defined.

\noindent Let $u_{n}$ be a sequence of $\mathcal{C}^{2}$--class functions
satisfying

$(i)$ \ $u_n$ converges uniformly to $u$ in the interval $[-R, R]$.

$(ii)$ \ $u_n^{\prime }$ converges uniformly to $u^{\prime }$ in the
interval $[-R, R]$

$(iii)$ \ $u_n^{\prime\prime}$ converges in $\mathbb{L}^{1}([-R, R])$ to $%
u^{\prime\prime}$.

\noindent We use It\^o's formula to show that,
\begin{equation*}
u_n(Y_{t\wedge\tau_{R}})=u_n(Y_{0})+\int_{0}^{t\wedge\tau_{R}}u_n^{\prime}
(Y_{s})dY_{s}+\frac{1}{2}\int_{0}^{t\wedge\tau_{R}}u_n^{\prime\prime}
(Y_{s})|Z_{s}|^{2}ds
\end{equation*}

%
%
%
%
%
%

\noindent Passing to the limit (on $n$) in the previous identity and using
the above properties $(i)$, $(ii)$, $(iii)$ and Proposition \ref{Krylovloc}
we get
\begin{equation*}
u(Y_{t\wedge\tau_{R}})=u(Y_{0})+\int_{0}^{t\wedge\tau_{R}}u^{\prime}
(Y_{s})dY_{s}+\frac{1}{2}\int_{0}^{t\wedge\tau_{R}}u^{\prime\prime}
(Y_{s})|Z_{s}|^{2}ds
\end{equation*}
Indeed, the limit for the left hand side term, as well as those of the first
and the second right hand side terms can be obtained by using properties $%
(i) $ and $(ii)$. The limit for the third right hand side term follows from
property $(iii)$ and Proposition \ref{Krylovloc}. \eop


\section{QBSDEs with $\mathbb{L}^2$ terminal data and measurable generators}

The present section will be developed in the same spirit of section 3. The It\^o--Krylov formula (established in section 4) will replace the It\^o formula in
all proofs. Thanks to It\^o--Krylov's formula, the following lemma, will play
the same role as Lemma \ref{zvonkin0C} when $f$ is merely measurable. In
particular, it allows us to eliminate the additive quadratic term from the
simple QBSDEs $eq(\xi ,f(y)|z|^{2})$ and $eq(\xi ,a+b|y|+c|z|+f(y)|z|^{2})$.

\begin{lemma}
\label{zvonkin0} Let $f$ belongs to $\mathbb{L}^{1}(\mathbb{R})$. The
function
\begin{equation}
u(x):=\int_{0}^{x}\exp \left( 2\int_{0}^{y}f(t)dt\right) dy
\end{equation}
satisfies then the following properties,

\noindent (i) \ $u\in \mathcal{C}^{1}(\mathbb{R})\cap \mathcal{W}_{1,loc}^{2}(\mathbb{R%
})$ and satisfies, for $a.e.\ x$, $\frac{1}{2}u^{\prime \prime
}(x)-f(x)u^{\prime }(x)=0$.

\noindent (ii) \ $u$ is a one to one function from $\mathbb{R}$ onto $\mathbb{R}$.

\noindent (iii) The inverse function $u^{-1}$ belongs to $\mathcal{C}^{1}(\mathbb{R}%
)\cap \mathcal{W}_{1,loc}^{2}(\mathbb{R})$.

\noindent (iv) \ There exist two positive constants $m$ and $M$ such that,

\hskip 1.5cm for any \ $x,y\in \mathbb{R}$, \ $m\left\vert x-y\right\vert
\leq \left\vert u(x)-u(y)\right\vert \leq M\left\vert x-y\right\vert $

\end{lemma}



\subsection{The equation $eq(\protect\xi, f(y)|z|^2)$}

\vskip0.15cm The following proposition shows that neither the exponential
moment of $\xi$ nor the continuity of the generator are needed to obtain the
existence and uniqueness of the solution to quadratic BSDEs.

\begin{proposition}
\label{f(y)|z|^2} Assume \textbf{\ (H1)} be satisfied. Let $f$ be a globally
integrable function on $\mathbb{R}$. Then, the BSDE $eq(\xi,f(y)|z|^{2})$
has a unique solution $(Y, \, Z)$ which belongs to $\mathcal{S}^2\times%
\mathcal{M}^2$.
\end{proposition}

\bop Let $u$ be the function defined in Lemma \ref{zvonkin0}. Since $\xi $
is square integrable then $u(\xi )$ is square integrable too. Therefore $%
eq(u(\xi ),\,0)$ has a unique solution in $\mathcal{S}^2\times\mathcal{M}^2$%
. The proposition follows now by applying It\^o--Krylov's formula to the
function $u^{-1}$ which belongs to the space $\mathcal{C}^{1}(\mathbb{R}
)\cap \mathcal{W}_{1,loc}^{2}(\mathbb{R})$. \eop

\vskip0.2cm The following proposition allows us to compare the solutions for
QBSDEs of type $eq(\xi ,f(y)|z|^{2})$. The novelty is that the comparison
holds whenever we can only compare the generators for $a.e.$ $y$. Moreover,
both the generators can be non--Lipschitz.

\begin{proposition}
\label{ct}(Comparison) Let $\xi_1$, $\xi_2$ be $\mathcal{F}_T$--measurable
and satisfy assumption \textbf{(H1)}. Let $f$, $g$ be in $\mathbb{L}^{1}(%
\mathbb{R})$. Let $( Y^{f},Z^{f}) $, $( Y^{g},Z^{g}) $ be respectively the
solution of the BSDEs $eq(\xi_1, f(y)|z|^2)$ and $eq(\xi_2, g(y)|z|^2)$.
Assume that $\xi_{1}\leq\xi_{2}$ a.s. and $f \leq g$ $a.e$. Then $%
Y_{t}^{f}\leq Y_{t}^{g}$ for all $t$ $\mathbb{P}$--a.s.
\end{proposition}

\bop According to Proposition \ref{f(y)|z|^2}, the solutions $( Y^{f},Z^{f})
$ and $( Y^{g},Z^{g}) $ belong to $\mathcal{S}^2\times\mathcal{M}^2$. For a
given function $h$, we put
\begin{equation*}
u_{h}(x):=\int_{0}^{x}\exp\left ( 2\int_{0}^{y}h(t)dt\right) dy
\end{equation*}
The idea consists to apply suitably Proposition \ref{Ito-Krylov} to the $%
u_{f}(Y_{T}^{g})$, this gives
\begin{align*}
u_{f}(Y_{T}^{g})& =u_{f}(Y_{t}^{g})+\int_{t}^{T}u_{f}^{\prime
}(Y_{s}^{g})dY_{s}^{g}+\dfrac{1}{2}\int_{t}^{T}u_{f}^{\prime \prime
}(Y_{s}^{g})d\langle Y_{\cdot }^{g}\rangle _{s} \\
& =u_{f}(Y_{t}^{g})+M_{T}-M_{t}-\int_{t}^{T}u_{f}^{\prime
}(Y_{s}^{g})g(Y_{s}^{g})|Z_{s}^{g}|^{2}ds \\
& \;\; +\dfrac{1}{2}\int_{t}^{T}u_{f}^{\prime \prime
}(Y_{s}^{g})|Z_{s}^{g}|^{2}ds
\end{align*}%
Since $u_{g}^{\prime \prime }(x)-2g(x)u_{g}^{\prime }(x)=0$, $u_{f}^{\prime
\prime }(x)-2f(x)u_{f}^{\prime }(x)=0$ and $u_{f}^{\prime }(x)\geq 0$, then
\begin{equation*}
u_{f}(Y_{T}^{g})=u_{f}(Y_{t}^{g})+M_{T}-M_{t}-\int_{t}^{T}u_{f}^{\prime
}(Y_{s}^{g})\left[ g(Y_{s}^{g})-f(Y_{s}^{g})\right] |Z_{s}^{g}|^{2}ds
\end{equation*}%
where $(M_{t})_{0\leq t\leq T}$ is a martingale. \newline
Since the term
\begin{equation*}
\int_{t}^{T}u_{f}^{\prime }(Y_{s}^{g})\left[ g(Y_{s}^{g})-f(Y_{s}^{g})\right]
|Z_{s}^{g}|^{2}ds
\end{equation*}
is positive, then
\begin{align*}
u_{f}\left( Y_{t}^{g}\right) \geq u_{f}(Y_{T}^{g})+M_{T}-M_{t}
\end{align*}
Since $Y_{t}^{f}$ and $Y_{t}^{g}$ is $\mathcal{F}_{t}$--adapted, then
passing to conditional expectation and using the fact that $u_{f}$ is an
increasing function and $\xi _{2} \geq \xi _{1}$, we get
\begin{align*}
u_{f}\left( Y_{t}^{g}\right) & \geq \mathbb{E}\left[ u_{f}(Y_{T}^{g})\left/
\mathcal{F}_{t}\right. \right] \\
& =\mathbb{E}\left[ u_{f}(\xi _{2})\left/ \mathcal{F}_{t}\right. \right] \\
& \geq \mathbb{E}\left[ u_{f}(\xi _{1})\left/ \mathcal{F}_{t}\right. \right]
\\
& = u_{f}\left( Y_{t}^{f}\right)
\end{align*}
Passing to $u_{f}^{-1}$, we get $Y_{t}^{g}\geq Y_{t}^{f}$. Proposition \ref%
{ct} is proved. \eop


\vskip 0.15cm The following uniqueness result is a consequence of the
previous proposition.

\begin{corollary}
\label{cct} Let $\xi $ satisfies \textbf{(H1)} and $f$, $g$ be integrable
functions. Let $(Y^{f},Z^{f})$ and $(Y^{g},Z^{g})$ respectively denote the
(unique) solutions of the BSDE $eq(\xi ,f(y)|z|^2)$ and $eq(\xi ,g(y)|z|^2)$
. If $f=g$--a.e., then $(Y^{f},Z^{f})=(Y^{g},Z^{g})$ in $\mathcal{S}%
^{2}\times\mathcal{M}^{2}$.
\end{corollary}

\begin{remark}
Proposition \ref{ct} and Corollary \ref{cct} will be used in the PDEs part,
to show the existence of a gap in the classical relation between the BSDEs
and their corresponding PDEs.
\end{remark}


\subsection{The equation $eq(\protect\xi, a+b|y|+c|z|+f(y)|z|^{2})$}

The BSDE under consideration in this subsection is,
\begin{equation}
Y_{t}=\xi
+\int_{t}^{T}(a+b|Y_{s}|+c|Z_{s}|+f(Y_{s})|Z_{s}|^{2})ds-%
\int_{t}^{T}Z_{s}dW_{s}  \label{edsrab2}
\end{equation}%
where $a, \, \, b,\, \, c \in \mathbb{R}$ and $f:\mathbb{R}\longmapsto
\mathbb{R}$.

\vskip0.2cm We refer to BSDE (\ref{edsrab2}) as equation $eq(\xi
,a+b|y|+c|z|+f(y)|z|^{2})$.

\begin{proposition}
\label{equivabf} Assume that \textbf{\ (H1)} is satisfied. Assume moreover
that $f$ is globally integrable on $\mathbb{R}$. Then, the BSDE $eq(\xi
,\,a+b|y|+ c| z|+f(y)|z|^{2})$ has a minimal and a maximal solution.
Moreover all solutions are in $\mathcal{S}^{2}\times\mathcal{M}^{2}$.
\end{proposition}


\bop Let $u$ be the function defined in Lemma \ref{zvonkin0}. Consider the
BSDE
\begin{equation}  \label{edsrLSM2}
\bar Y_{t}=\bar \xi +\int_{t}^{T}G(\bar Y_{s}, \bar Z_{s})ds-
\int_{t}^{T}\bar Z_{s}dW_{s}
\end{equation}
where $G(\bar y,\bar z):=(a+b|u^{-1}(\bar y)|)u^{\prime }[u^{-1}(\bar y)] +
c|\bar z|$.

\noindent From Lemma \ref{zvonkin0}, we deduce that the terminal condition $%
\bar\xi :=u(\xi )$ is square integrable (since Assumption \textbf{(H1)}) and
the generator $G$ is continuous and with linear growth. Arguing then as in
the "alternative proof of Proposition \ref{equivabfC}", one can prove that
the BSDE (\ref{edsrLSM2}) has a maximal and a minimal solutions in $\mathcal{%
S}^{2}\times\mathcal{M}^{2} $. Applying now the It\^o--Krylov formula to the
function $u^{-1}(\bar Y_t)$, we show that the BSDE $eq(\xi ,\,a+b|y|+ c|
z|+f(y)|z|^{2})$ has a solution. Since $u$ is a strictly increasing
function, we then deduce the existence of a minimal and a maximal solutions
for the initial equation $eq(\xi ,\,a+b|y|+ c| z|+f(y)|z|^{2})$. \eop


\subsection{The BSDE $eq (\protect\xi ,H)$ with $\vert H(s,y,z)\vert \leq a
+ b\vert y\vert + c\vert z\vert + f(\vert y\vert)|z|^2$}

\noindent Consider the assumption,

\vskip 0.2cm\noindent \textbf{(H6)} There exist positive real numbers $a,b,
c $ such that for every $s,y,z$ \vskip 0.2cm \hskip 2cm $|H(s,y,z)|\leq
a+b|y|+ c\vert z\vert + f(|y|)|z|^{2}:=g(y,z)$, \vskip 0.2cm\noindent where $%
f$ is some positive locally bounded integrable function, but not necessarily
continuous.

\begin{theorem}
\label{BSDEHdominé} Assume that (\textbf{H1}), (\textbf{H4}) and (\textbf{H6}%
) are fulfilled. Then, the BSDE $eq(\xi ,H)$ has at least one solution \ $%
(Y,Z)$ which belongs to $\mathcal{S}^{2}\times\mathcal{M}^{2}$.
\end{theorem}



\begin{remark}
\label{flocalementborne} Although the proof of Theorem \ref{BSDEHdominé} may
be performed as that of Theorem \ref{BSDEHdominéC}, we will give another
proof which consists to use a comparison theorem and an appropriate
localization by a suitable dominating process which is derived from the
extremal solutions of the two QBSDEs $eq(-\xi^- ,
-(a+b|y|+c|z|+f(|y|)|z|^{2}))$ and $eq(\xi^+ ,(a+b|y|+c|z|+f(|y|)|z|^{2}))$.
\end{remark}

\vskip 0.15cm To prove Theorem \ref{BSDEHdominé}, we need the following two
lemmas. The first one allows us to show that $Z$ belongs to $\mathcal{M}^2$
while the second is a comparison theorem for our context.


\begin{lemma}
\label{zvonkin1} Let $f$ belongs to $\mathbb{L}^{1}(\mathbb{R})$ and put \ $%
K(y):=\int _{0}^{y} \exp(-2\int_{0}^{x}f(r)dr)dx$. The function $$u(x):=\int_{0}^{x}K(y) \exp\left(2\int_{0}^{y}f(t)dt\right)dy$$ satisfies following
properties:

\vskip 0.15cm \noindent $(i)$ $u$ belongs to $\mathcal{C}^{1}(\mathbb{R})\cap\mathcal{W%
} _{1,loc}^{2}(\mathbb{R})$, and, $u(x)\geq0$ and $u^{\prime }(x)\geq0$ for $%
x\geq0$.

\vskip 0.15cm \noindent Moreover, $u$ satisfies, for a.e. $x$, $\frac{1}{2}%
u^{\prime\prime}(x)- f(x)u^{\prime}(x)= \frac{1}{2}$.

\vskip 0.15cm \noindent $(ii)$ The map $x\mapsto v(x) := u(\vert x\vert)$ belongs
to $\mathcal{C} ^{1}(\mathbb{R})\cap\mathcal{W} _{1,loc}^{2}(\mathbb{R})$,
and $v^{\prime}(0) = 0$.

\vskip 0.15cm \noindent $(iii)$ There exist $c_1>0,\ c_2>0$ such that for every $x\in%
\mathbb{R}$, $u(\vert x\vert) \leq c_1 \vert x\vert^2$ and $u^{\prime}(\vert
x\vert) \leq c_2\vert x\vert $.
\end{lemma}

\begin{lemma}
\label{comparison} (Comparison) Let $h_1(t,\omega, y, z)$ be uniformly
Lipschitz in $(y,z)$ uniformly with respect to $(t,\omega)$. Let $%
h_2(t,\omega, y, z)$ be $\mathcal{F}_t$--progressively measurable and such
that for every process $(U,V)\in\mathcal{S}^2\times\mathcal{M}^2$, $\mathbb{E%
}\int_0^T \vert h_2(s,U_s,V_s)\vert ds<\infty$. Let $\xi _{i}\in
L^{2}(\Omega )$, $(i=1,2)$, be an $\mathcal{F}_{T}$--measurable random
variables. Let $(Y^{1},\ Z^{1})\in \mathcal{S}^{2}\times \mathcal{M}^{2}$ be
the unique solution of BSDE $eq(\xi _{1},h_{1})$ and $(Y^{2},\,Z^{2})\in
\mathcal{S}^{2}\times \mathcal{M}^{2}$ be a solution of BSDE $eq(\xi
_{2},h_{2})$. Assume that, $\xi _{1}\leq \xi _{2}$ for $a.s.\ \omega $ and $%
h_{1}(s,Y_{s}^{2},Z_{s}^{2})\leq h_{2}(s,Y_{s}^{2},Z_{s}^{2})$ for $a.e. $ $%
s,\omega $. Then, $Y_{t}^{1}\leq Y_{t}^{2}$ \ for every $t$ and $a.s.\
\omega $.
\end{lemma}

\bop Applying It\^o's formula to $(( Y_{t}^{1}-Y_{t}^{2}) ^{+}) ^{2}$, and since $(\xi _{1}-\xi _{2})^{+}=0$,
\begin{align*}
\left( \left( Y_{t}^{1}-Y_{t}^{2}\right) ^{+}\right) ^{2}& +\int_{t}^{T}%
\mathbf{1}_{\left\{ Y_{s}^{1}>Y_{s}^{2}\right\} }\left\vert
Z_{s}^{1}-Z_{s}^{2}\right\vert ^{2}ds \\
& =2\int_{t}^{T}\left( Y_{s}^{1}-Y_{s}^{2}\right) ^{+}\left[ h_{1}\left(
s,Y_{s}^{1},Z_{s}^{1}\right) -h_{2}\left( s,Y_{s}^{2},Z_{s}^{2}\right) %
\right] ds \\
& \ -2\int_{t}^{T}\left( Y_{s}^{1}-Y_{s}^{2}\right) ^{+}\left[
Z_{s}^{1}-Z_{s}^{2}\right] dW_{s}
\end{align*}%
Passing to expectation and using the fact that $h_{1}(s,Y_{s}^{2},Z_{s}^{2})%
\leq h_{2}(s,Y_{s}^{2},Z_{s}^{2})$, we obtain
\begin{align*}
& \mathbb{E}\left[ \left( Y_{t}^{1}-Y_{t}^{2}\right) ^{+}\right] ^{2} +%
\mathbb{E}\int_{t}^{T}\mathbf{1}_{\left\{ Y_{s}^{1}>Y_{s}^{2}\right\}
}\left\vert Z_{s}^{1}-Z_{s}^{2}\right\vert ^{2}ds \\
& =2\mathbb{E}\int_{t}^{T}\left( Y_{s}^{1}-Y_{s}^{2}\right) ^{+}\left[
h_{1}\left( s,Y_{s}^{1},Z_{s}^{1}\right) -h_{1}\left(
s,Y_{s}^{2},Z_{s}^{2}\right) \right] ds \\
& \quad +2\mathbb{E}\int_{t}^{T}\left( Y_{s}^{1}-Y_{s}^{2}\right) ^{+}\left[
h_{1}\left( s,Y_{s}^{2},Z_{s}^{2}\right) -h_{2}\left(
s,Y_{s}^{2},Z_{s}^{2}\right) \right] ds \\
& \leq 2\mathbb{E}\int_{t}^{T}\left( Y_{s}^{1}-Y_{s}^{2}\right) ^{+}\left[
h_{1}\left( s,Y_{s}^{1},Z_{s}^{1}\right) -h_{1}\left(
s,Y_{s}^{2},Z_{s}^{2}\right) \right] ds \\
& \leq 2\mathbb{E}\int_{t}^{T}\left( Y_{s}^{1}-Y_{s}^{2}\right) ^{+}[\mathbf{%
1}_{\left\{ Y_{s}^{1}>Y_{s}^{2}\right\} }+\mathbf{1}_{\left\{ Y_{s}^{1}\leq
Y_{s}^{2}\right\} }]\\ & \qquad \qquad \times \left[ h_{1}\left( s,Y_{s}^{1},Z_{s}^{1}\right)
-h_{1}\left( s,Y_{s}^{2},Z_{s}^{2}\right) \right] ds
\end{align*}%
Since $\mathbf{1}_{\left\{ Y_{s}^{1}\leq Y_{s}^{2}\right\} }\left(
Y_{s}^{1}-Y_{s}^{2}\right) ^{+}=0$, we then have
\begin{align*}
& \mathbb{E}\left[ \left( Y_{t}^{1}-Y_{t}^{2}\right) ^{+}\right] ^{2}+ \
\mathbb{E}\int_{t}^{T}\mathbf{1}_{\left\{ Y_{s}^{1}>Y_{s}^{2}\right\}
}\left\vert Z_{s}^{1}-Z_{s}^{2}\right\vert ^{2}ds \\
& \leq 2\mathbb{E}\int_{t}^{T}\left( Y_{s}^{1}-Y_{s}^{2}\right) ^{+}\mathbf{1%
}_{\left\{ Y_{s}^{1}>Y_{s}^{2}\right\} }\left[ h_{1}\left(
s,Y_{s}^{1},Z_{s}^{1}\right) -h_{1}\left( s,Y_{s}^{2},Z_{s}^{2}\right) %
\right] ds
\end{align*}%
Using the fact that $h_{1}$ is Lipschitz and the inequality $ab\leq
\varepsilon a^{2}+\frac{1}{\varepsilon }b^{2}$, we obtain
\begin{equation*}
E\left\vert \left( Y_{t}^{1}-Y_{t}^{2}\right) ^{+}\right\vert ^{2}\leq
CE\int_{t}^{T}\left\vert (Y_{s}^{1}-Y_{s}^{2})^{+}\right\vert ^{2}ds.
\end{equation*}
Using Gronwall's lemma, we get \ $\left( Y_{t}^{1}-Y_{t}^{2}\right) ^{+}=0$
for every $t$ and $a.s.$ $\omega$, which implies that $Y_{t}^{1}\leq
Y_{t}^{2}$ \ for every $t$, and $a.s.$ $\omega$. \eop

\vskip0.2cm \noindent \textbf{Proof of Theorem \ref{BSDEHdominé}} We assume
for simplicity that $\xi $ is positive. Let $Y^{g}$ be the maximal solution
of the BSDE $eq(\xi ,g)$ and $Y^{-g}$ be the minimal solution of the BSDE $%
eq(-\xi ,-g)$. By Proposition \ref{equivabf}, $Y^{-g}$ and $Y^{g}$ exist and
belong to $\mathcal{S}^{2}\times \mathcal{M}^{2}$. Let $\xi _{n}:=\xi \wedge
n$. Let $(H_{n})$ be an increasing sequence of Lipschitz functions which
converges to $H$ uniformly on compact sets. For each $n$, we denote by $%
(Y^{n},Z^{n})$ the unique solution of the BSDE $eq(\xi _{n},H_{n})$. We know
that for every $n$, $(Y^{n},Z^{n})$ belongs to $\mathcal{S}^{2}\times
\mathcal{M}^{2}$. Since \ $\xi _{n}\leq \xi $ and $H_{n}\leq g$ \ for each $n
$, then Lemma \ref{comparison} (comparison) shows that for every $n$, $t$
and $a.s.\ \omega $,
\begin{equation}
|Y_{t}^{n}|\leq |Y_{t}^{-g}|+|Y_{t}^{g}|:=S_{t}  \label{S}
\end{equation}%
For $R>0$, we define a stopping time $\tau _{R}$ by
\begin{equation}
\tau _{R}:=\inf \{t\geq 0\ :\ S_{t}\geq R\}\wedge T  \label{tauR}
\end{equation}%
The process \ $(Y_{t}^{n,R},\,Z_{t}^{n,R}):=(Y_{t\wedge \tau _{R}}^{n},%
\mathbf{1}_{\{t\leq \tau _{R}\}}Z_{t}^{n})_{0\leq t\leq T}$ \ satisfies then
the BSDE $eq(Y_{\tau _{R}}^{n,R},H_{n},R)$
\begin{eqnarray*}
Y_{t}^{n,R}&=&Y_{\tau _{R}}^{n,R}+\int_{t}^{T}\mathbf{1}_{\{s\leq \tau
_{R}\}}H_{n}(s,Y_{s}^{n,R},Z_{s}^{n,R})ds -\int_{t}^{T}Z_{s}^{n,R}dW_{s}\notag
 \end{eqnarray*}
From inequality (\ref{S}) and the definition of $\tau _{R}$, we deduce that
for every $n$ and every $t\in \lbrack 0,\ T]$,
\begin{equation}
|Y_{t}^{n,R}|\ \leq \ S_{t\wedge \tau _{R}}\leq R  \label{borneYn}
\end{equation}%
As in \cite{BHM} (see also \cite{LSM2}), we define a function $\rho $ by,
\begin{equation*}
\rho (y):=-R\mathbf{1}_{\{y<-R\}}+y\mathbf{1}_{\{-R\leq y\leq R\}}+R\mathbf{1%
}_{\{y>R\}}
\end{equation*}%
\noindent It is not difficult to prove that $(Y^{n,R},\,Z^{n,R})$ solves the
BSDE \\ $eq(Y_{\tau _{R}}^{n,R},H_{n}(t,\rho (y),z))$. \newline
Since $f$ is locally bounded, then for every $y$ satisfying $|y|\leq R$ we
have,
\begin{align}
|H_{n}(t,\rho (y),z)|& \leq a+b|\rho (y)|+c|z|+\sup_{|y|\leq R}f(|\rho
(y)|)|z|^{2}  \notag  \label{bornHn} \\
& \leq a_{1}+bR+C(R)|z|^{2}
\end{align}%
where $a_{1}:=a+c^{2}$ \ and \ $\displaystyle C(R)=1+\sup_{|y|\leq R}f(|\rho
(y)|)$. \newline
Therefore, passing to the limit on $n$ and using the Kobylanski monotone
stability result \cite{Kob}, one can show that for any $R>0$, the sequence $%
(Y^{n,R},Z^{n,R})$ converges to a process $(Y^{R},Z^{R})$ which satisfies
the following BSDE on $[0, \ T]$,
\begin{equation}
Y_{t}^{R}=\xi ^{R}+\int_{t}^{T}\mathbf{1}_{\{s\leq \tau
_{R}\}}H(s,Y_{s}^{R},Z_{s}^{R})ds-\int_{t}^{T}Z_{s}^{R}dW_{s},
\tag{$eq(\xi
,H,R)$}  \label{edsrYR}
\end{equation}%
where $Y_{t}^{R}:= \sup_{n}Y_{t}^{n,R}$ and $\xi ^{R}:= \sup_{n}Y_{\tau
_{R}}^{n}$.

\vskip 0.15cm Arguing as in the proof of Theorem \ref{BSDEHdominéC}, but we
use Lemma \ref{zvonkin1} (in place of Lemma \ref{zvonkin1C}) and It\^o%
--Krylov's formula (in place of It\^o's formula), one can show that \ $%
(Y^{R},\,Z^{R})$ \ belongs to $\mathcal{S}^{2}\times \mathcal{M}^{2}$. \vskip%
0.15cm \noindent We define
\begin{equation*}
Y_{t}:=\lim_{R\rightarrow \infty }Y_{t\wedge \tau _{R}}\ \ \text{for}\ t\in
\lbrack 0,T]
\end{equation*}%
and
\begin{equation*}
Z_{t}:=Z_{t}^{R}\ \ \text{for}\ t\in (0,\tau _{R})
\end{equation*}

\noindent Passing to the limit on $R$ in the previous BSDE $eq(\xi, H, R)$,
one can show that the process $(Y,Z):= (Y_t,Z_t)_{t\leq T}$ satisfies the
BSDE $eq(\xi, H)$.

Since $Y^{-g}$ and $Y^{g}$ belong to $\mathcal{S}^2$, we deduce that $Y$
belongs to $\mathcal{S}^2$ also.

Arguing as in the proof of Theorem \ref{BSDEHdominéC}, but use Lemma \ref%
{zvonkin1} (in place of Lemma \ref{zvonkin1C}) and It\^o--Krylov's formula (in
place of It\^o's formula) one can prove that $Z\in \mathcal{M}^{2}$.

We shall prove that $Y_{T} = \xi$. Let $Y^{\prime}$ be the minimal solution
of BSDE $eq(\xi, -g)$ and $Y^{\prime\prime}$ be the maximal solution of BSDE
$eq(\xi, g)$. Using Lemma \ref{comparison}, we obtain, for any $n$, $R$ and $%
t\in [0, \ T]$
\begin{equation*}
Y_{t\wedge\tau_R}^{\prime} \leq Y_{t\wedge\tau_R}^{n} \leq
Y_{t\wedge\tau_R}^{\prime\prime}
\end{equation*}
Passing to the limit on $n$, we obtain,
\begin{equation*}
Y_{t\wedge\tau_R}^{\prime} \leq \sup_{n} Y_{t\wedge\tau_R}^{n} :=
Y_{t\wedge\tau_R} \leq Y_{t\wedge\tau_R}^{\prime\prime}
\end{equation*}
Putting $t=T$, we get for any $R$
\begin{equation*}
Y_{\tau_R}^{\prime} \leq Y_{\tau_R} \leq Y_{\tau_R}^{\prime\prime}
\end{equation*}
Since $Y^{\prime}$ and $Y^{\prime\prime}$ are continuous and $Y_{T}^{\prime}
= Y_{T}^{\prime\prime} = \xi$, then letting $R$ tends to infinity, we get
\begin{equation*}
\xi \leq \liminf_{R\rightarrow\infty} Y_{\tau_R} \leq
\limsup_{R\rightarrow\infty} Y_{\tau_R} \leq \xi
\end{equation*}
Theorem \ref{BSDEHdominé} is proved. \eop

\begin{remark}
We are currently work to drop the global integrability condition on the
function $f$. The situation becomes more delicate under this (local
integrability) condition. It requires some localization arguments and
supplementary assumptions on the terminal condition.
\end{remark}


\section{Application to Quadratic Partial Differential Equations} \ \\
Let $\sigma $, $b$ be measurable functions defined on $\mathbb{R}^{d}$ with
values in $\mathbb{R}^{d\times d}$ and $\mathbb{R}^{d}$ respectively.

Let $a:=\sigma \sigma ^{\ast }$ and define the operator $L$ by
\begin{equation*}
{L}:=\sum_{i,\,j=1}^{d}{a}_{ij}(x)\frac{\partial ^{2}}{\partial
x_{i}\partial x_{j}}+\sum_{i=1}^{d}{b}_{i}(x)\frac{\partial }{\partial x_{i}}
\end{equation*}%
Let $\psi $ be a measurable function from $\mathbb{R}^{d}$ to $\mathbb{R}$.
Consider the following semi-linear PDE
\begin{equation}
\left\{
\begin{array}{l}
\dfrac{\partial {v}}{\partial s}(s,\,x)={L}v(s,x)+f(v(s,\,x))|\nabla
_{x}v(s,\,x))|^{2},\ \text{on\ } [0,T)\times \mathbb{R}^{d} \\
\\
v(T,x)=\psi (x)%
\end{array}%
\right.  \label{edpinitiale}
\end{equation}

\noindent {\textbf{Assumptions.}}

\textbf{(H7)} $\sigma$, $b$ are uniformly Lipschitz.


\textbf{(H8)} $\sigma$, $b$ are of linear growth and $f$ is continuous and
integrable.

\textbf{(H9)} The terminal condition $\psi$ is continuous and with
polynomial growth.

\begin{theorem}
\label{viscosity} Assume \textbf{(H7), (H8)} and (\textbf{H9}) hold. Then, $%
v(t,x):=Y_{t}^{t,x}$ is a viscosity solution for the PDE (\ref{edpinitiale}).
\end{theorem}

\begin{remark}
The conclusion of Theorem \ref{viscosity} remains valid if the assumption
\textbf{(H7)} is replaced by:``the martingale problem is well--posed for $%
a:=\sigma\sigma^*$ and $b$".
\end{remark}

To prove the existence of viscosity solution, we will follow the idea of
\cite{Kob}. To this end, we need the following touching property. This
allows to avoid the comparison theorem. The proof of the touching property
can be found for instance in \cite{Kob}.

\begin{lemma}
Let $(\xi _{t})_{0\leq t\leq T}$ be a continuous adapted process such that
\begin{equation*}
d\xi _{t}=\beta (t)dt+\alpha (t)dW_{t},
\end{equation*}
where $\beta $ and $\alpha $ are continuous adapted processes such that $b$,
$|\sigma |^{2}$ are integrable. If $\xi _{t}\geq 0$ a.s. for all $t$, then
for all $t$,
\begin{equation*}
\mathbf{1}_{\{\xi _{t}=0\}}\alpha (t)=0\ \ \ a.s.,
\end{equation*}
\begin{equation*}
\mathbf{1}_{\{\xi _{t}=0\}}\beta (t)\geq 0\ \ \ a.s.,
\end{equation*}
\end{lemma}

\vskip 0.2cm\noindent \textbf{Proof of Theorem \ref{viscosity}.} \ We first
prove the continuity of $v(t,x):=Y_t^{t,x}$. Let $u$ be the transformation
defined in Lemma \ref{zvonkin0}. Let $(\bar Y_s^{t,x}, \bar Z_s^{t,x})$ be
the unique $\mathcal{M}^2$ solution of the BSDE $eq(u(\psi(X_{T}^{t,x}), 0)$%
. Using assumption \textbf{(H7)}, one can show that the map $(t,x)
\mapsto \bar Y_t^{t,x}$. Using Lemma \ref{zvonkin0}, we deduce that $%
v(t,x):=Y_t^{t,x}$ is continuous in $(t,x)$. We now show that $v$ is a
viscosity subsolution for PDE (\ref{edpinitiale}). We denote $(X_{s},
Y_{s},Z_{s }) := (X_{s}^{t,x}, Y_{s}^{t,x},Z_{s}^{t,x})$. Since $%
v(t,x)=Y_{t}^{t,x}$, then the Markov property of $X$ and the uniqueness of $%
Y $ show that for every $s\in [0,T]$
\begin{equation}
v(s,X_{s})=Y_{s}  \label{markov}
\end{equation}

Let $\phi \in \mathcal{C}^{1,2}$ and $(t,x)$ be a local maximum of $v-\phi $
which we suppose global and equal to $0$, that is :
\begin{equation*}
\phi (t,x)=v(t,x) \text{ \ \ \ and \ \ \ } \phi (\overline{t},\overline{x }%
)\geq v(\overline{t},\overline{x}) \text{ \ for each \ } (\overline{t},%
\overline{x}).
\end{equation*}

\vskip0.1cm\noindent This and equality (\ref{markov}) imply that
\begin{equation}
\phi (s,X_{s})\geq Y_{s}  \label{phi>Y}
\end{equation}
By It\^o's formula we have
\begin{equation*}
\phi (s,X_{s})=\phi (t,X_{t})+\int_{t}^{s}\left( \frac{\partial \phi }{%
\partial r}+L\phi \right) (r,X_{r})dr+\int_{t}^{s}\sigma \nabla_x\phi
(r,X_{r})dW_{r}
\end{equation*}
Since $\phi (s,X_{s})\geq Y_{s}$, and $Y$ satisfies the equation
\begin{equation*}
Y_{t}=Y_{s}+\int_{t}^{s}f(Y_{r})\vert Z_{r}\vert^2
dr-\int_{t}^{s}Z_{r}dW_{r},
\end{equation*}
then the touching property shows that for each $s $,
\begin{equation*}
\mathbf{1}_{\{\phi (s,X_{s})=Y_{s}\}}\left( \frac{\partial \phi }{\partial t}
+L\phi \right) (s,X_{s})+f(Y_{s})\vert Z_{s}\vert^2\geq 0\ \ \ a.s.,
\end{equation*}
and
\begin{equation*}
\mathbf{1}_{\{\phi (s,X_{s})=Y_{s}\}}|\sigma ^{T}\nabla_x\phi
(s,X_{s})-Z_{s}|\geq 0\ \ \ a.s.
\end{equation*}

\noindent Since for $s=t$, $\phi (t,x):=\phi (t,X_{t})=Y_{t}:= v(t,x)$, then
the second equation gives $Z_{t}=\sigma \nabla_x\phi (t,X_{t}):=\sigma
\nabla_x\phi (t,x)$, and the first inequality gives the desired result.
\eop

\begin{remark} \label{gap}
From Theorem \ref{viscosity}, one can see that there is a gap between the
solution of a BSDE and the viscosity solutions of its associated PDE. That
is, the existence of a unique solution to a BSDE (even when the comparison
theorem holds) does not systematically allow to define a viscosity solution
to the associated PDE. Indeed, the Corollary \ref{cct} shows that the QBSDE $%
eq(\xi ,f_{1}(y)|z|^{2})$ and $eq(\xi ,f_{2}(y)|z|^{2})$ generate the same
solution when $f_{1}$ and $f_{2}$ are equal almost surely. Thereby, for a
square integrable $\xi $ and $f$ belonging to $\mathbb{L}^{1}(\mathbb{R})$,
the QBSDE $eq(\xi ,f(y)|z|^{2})$ has a unique solution in $\mathcal{S}%
^{2}\times \mathcal{M}^{2}$, but how define the associated PDE (\ref%
{edpinitiale}) when $f$ is defined merely $a.e.$ ? What meaning to give to $%
f(v(t,x))$ when $v(t,x)$ stays to the set where $f$ is not defined ? We
think that, when $f$ is defined merely $a.e.$, the associated PDE associated
to BSDE $eq(\xi ,f_{2}(y)|z|^{2})$ would has the form
\begin{equation}
\left\{
\begin{array}{l}
\dfrac{\partial {v}}{\partial s}(s,\,x)={L}v(s,x)+f(v(s,\,x))|\nabla
_{x}v(s,\,x))|^{2},\ \text{on \ } [0,T)\times \mathbb{R}^{d} \\
\\
v(T,x)=\psi (x) \\
\\
\nabla v(t,x)=0\ \ \ if\ \ \ v(t,x)\in \mathcal{N}_{f}%
\end{array}%
\right.  \label{newedp}
\end{equation}%
where $\mathcal{N}_{f}$ denotes the negligible set of all real numbers $y$
for which $f$ is not defined.
\end{remark}

\end{document}